\documentclass[10pt]{article}
\usepackage[left=3cm,top=3cm,right=3cm,bottom=3cm]{geometry}
\usepackage{amsfonts,amssymb,amsmath,latexsym,indentfirst}
\usepackage{fontenc,enumerate}
\usepackage{graphics}
\usepackage{color,graphicx}
\usepackage[pdftex]{hyperref}
\usepackage{latexsym,wasysym,mathrsfs,enumerate}

\definecolor{darkred}{rgb}{.0,.0,.8}
\definecolor{darkgreen}{rgb}{.0,0,0.8}
\definecolor{darkredd}{rgb}{0,0,.8}
\hypersetup{
 colorlinks,%
       citecolor=darkgreen,%
    filecolor=magenta,%
    linkcolor=darkredd,%
    urlcolor=darkred,%
}
\numberwithin{equation}{section}
\newtheorem{proposition}{Proposition}[section]
\newtheorem{remark}{Remark}[section]
\newtheorem{definition}{Definition}[section]
\newtheorem{lemma}{Lemma}[section]
\newtheorem{theorem}{Theorem}[section]

\newcommand{\N}{\ensuremath{{\mathbb{N}} }}

\newcommand{\<}{\langle}
\renewcommand{\>}{\rangle}
\newcommand{\ga}{\gamma}
\newcommand{\eps}{\epsilon}

\newcommand{\si}{\sigma}

\newcommand{\al}{\alpha}

\newcommand{\rr}{{\mathbb R}}

\newcommand{\lam}{{\lambda}}
\newcommand{\ff}{\varphi}
\newcommand{\ti}{\widetilde}
\newcommand{\what}{\widehat}
\newcommand{\dd}{{\rm d}}
\newcommand{\ee}{{\rm e}}
\newcommand{\ii}{{\rm i}}
\newcommand{\h}{\mathscr{H}}

\newcommand{\kb}{\mathbb{K}}
\newcommand{\g}{\mathscr{G}(\beta,c,\gamma)}
\newcommand{\x}{\mathscr{X}}
\newcommand{\n}{\mathscr{N}}
\newcommand{\ddd}{\mathscr{D}}
\newcommand{\p}{\mathscr{P}}

\newcommand{\uu}{\mathscr{U}}
\newcommand{\dx}{D_x^{1/2}}
\newcommand{\st}{{H^{1/2}(\mathbb{R})}}
\newcommand{\dn}{\partial_x^{-1}}
\newcommand{\lt}{{L^2(\mathbb{R})}}

\newcommand{\fim}{\hfill$\square$\\ \\}
\newcommand{\e}{\varepsilon}
\newcommand{\q}{\quad}
\newcommand{\proof}{\noindent\textbf{Proof.}\quad}

\author{{\bf Amin Esfahani}\\  {\small School of Mathematics and Computer Science} \\
{\small Damghan University}\\ {\small Damghan, Postal Code 36716-41167, Iran}\\
{\small  E-mail: amin@impa.br, saesfahani@du.ac.ir}\vspace{2mm}\\
{\bf Steven Levandosky}\\
{\small Mathematics and Computer Science Department}\\
{\small College of the Holy Cross, Worcester, MA 01610} \\
{\small  E-mail: spl@mathcs.holycross.edu}}

\title{Solitary waves of the rotation-generalized Benjamin-Ono equation
}

\date{}

\begin{document}
\maketitle

\begin{abstract}
This work studies the rotation-generalized Benjamin-Ono equation which is derived from the theory of weakly nonlinear long surface and internal waves in deep water under the presence of rotation. It is shown that the solitary-wave solutions are orbitally stable for certain wave speeds.
\end{abstract}

\section{Introduction}
In the present paper we are concerned with studying the rotation-generalized Benjamin-Ono (RGBO) equation which can be written as
\begin{equation}\label{rgbo}
(u_t+\beta\h u_{xx}+(f(u))_x)_x=\gamma u,\q x\in\rr,\;t>0,
\end{equation}
where  $\ga>0$ and $\beta\neq0$ are real constants, $f$ is a $C^2$ function which is homogeneous of degree $p>1$ such that $sf(s)=pf'(s)$, and $\h$ denotes the Hilbert transform defined by
\[
\h u(x,t)={\rm{p.v.}}\frac{1}{\pi}\int_\rr\frac{u(z,t)}{x-z}\dd z,
\]
where p.v. denotes the Cauchy principal value. When $f(u)=\frac{1}{2}u^2$, equation \eqref{rgbo},  which is so-called the rotation-modified Benjamin-Ono (RMBO) equation, models the propagation of long internal waves in a deep rotating fluid \cite{galstep,grim,linmil,red}. In the context of shallow water the propagation of long waves in rotating fluid is described by the Ostrovsky equation \cite{belinov,ggs,ostrovsky,oststep}
\begin{equation}
(u_t+\beta u_{xxx}+(u^2)_x)_x=\gamma u,\q x\in\rr,\;t>0,
\end{equation}
which is also called the rotation-modified Korteweg-de Vries (RMKdV) equation. See also \cite{chen,esf} for the two-dimensional long internal waves in a rotating fluid. The parameter $\gamma$ is a measure of the effect of rotation. Setting $\gamma=0$ in \eqref{rgbo}, integrating the result with respect to $x$ and setting the constant of integration to zero, one obtains the generalized Benjamin-Ono (GBO) equation
\begin{equation}\label{gbo}
u_t+\beta\h u_{xx}+(f(u))_x=0.
\end{equation}

Most attention in this work is paid to the existence, the stability and the properties of localized traveling waves (commonly referred to as solitary waves) of \eqref{rgbo}.
Using variational methods and the Pohozaev-type identities, we prove the existence and nonexistence of solitary waves for a range of the parameters of \eqref{rgbo}. We also show that our solitary waves (of \eqref{rgbo}) are the ground states, i.e. they have minimal action. We also consider the effect of letting the rotation parameter $\ga$ approach zero. Actually we show that the ground state solitary waves converge to solitary waves of the GBO equation.

It was shown by Linares and Milanes \cite{linmil} that the RMBO equation \eqref{rgbo} is well-posed in the space
\[
X_s=\{f\in H^s(\rr); \partial_x^{-1}f\in H^s(\rr)\}
\]
with norm
\[
{\|g\|}_{X_s}=\|g\|_{H^s(\rr)}+\|\partial_x^{-1}g\|_{H^s(\rr)},
\]
for $s>3/2$, where the operator $\partial_x^{-1}$ is defined via the Fourier transform as $\what{\partial_x^{-1}g}(\xi)=(\ii\xi)^{-1}\what{g}(\xi)$. The methods therein also imply the same result for the RGBO equation \eqref{rgbo}.

It is also standard to show that the solution $u(t)$ obtained that way satisfies $E(u(t)) = E(u(0))$,
$Q(u(t)) = Q(u(0))$ and $M(u(t)) = 0$, for $t\in[0, T )$ with the maximum existence time $T$, where
\begin{equation}\label{cons-1}
E(u)=\int_\rr\frac{\beta}{2}(D^{1/2}_xu)^2+\frac{\ga}{2}(\partial_x^{-1}u)^2+F(u)\dd x,
\end{equation}
\begin{equation}\label{cons-2}
Q(u)=\frac{1}{2}\int_\rr u^2\dd x
\end{equation}
and
\begin{equation}
M(u)=\int_\rr u\;\dd x
\end{equation}
express, respectively, the energy, momentum and total mass, where $\what{\dx f}(\xi)=|\xi|^{1/2}\what{f}(\xi)$, $F'=f$ and $F(0)=0$. It is also worth remarking that the sufficiently smooth solutions of \eqref{rgbo} satisfy
\[
\int_\rr xu\;\dd x=0.
\]
These conserved quantities play an important role in our stability analysis.

We show in Theorems \ref{stability-theo} and \ref{instability-theo-2} that the function $d(c)$ defined by \eqref{d-function} determines the stability  of the solitary waves in the sense that if $d''(c)>0$, then $\g$ is $\x$-stable, while if $d''(c)<0$, then $\mathcal{O}_\varphi$ is $\x$-unstable, where the space $\x$ is defined in \eqref{x-space}. In Theorem \ref{main-instability}, we use the ideas of \cite{ribeiro}, and provides sufficient conditions for instability directly in terms of the parameters $\beta$, $\ga$ and $p$.

We also investigate the properties of the function $d(c)$ which determines the stability of the ground states. Using an important scaling identity, together with numerical approximations of the solitary waves, we are able to numerically approximate $d(c)$.

\begin{remark}
It is noteworthy that despite our regularity assumption on $f$, one can observe that all our results are valid for the nonlinearity $f(u)=-|u|u$.
\end{remark}
\section*{Notations}
For each $r\in\rr$, we define the translation operator by $\tau_ru= u(\cdot+r)$.

Given a solitary wave $\ff$ of \eqref{rgbo-1},
the orbit of $\ff$ is defined by the set $\mathcal{O}_\ff = \{\tau_r\ff;\;r\in\rr\}$.

We shall denote by $\widehat\ff$ the Fourier transform of $\ff$, defined as
\[
\widehat\ff(\zeta)=\int_{\rr}\;\ff(\omega)\ee^{-\ii \omega\cdot\zeta}\;\dd\omega.
\]
For $s\in\rr$, we denote by $H^s(\rr)$, the nonhomogeneous Sobolev space defined by the closure
\[
\left\{\ff\in\mathscr{S}'\left(\rr\right)\;:\;\|\ff\|_{H^s\left(\rr\right)}<\infty\right\},
\]
with respect to the norm
\[
\|\ff\|_{H^s\left(\rr\right)}
=\left\|\left(1+|\zeta|\right)^\frac{s}{2}\widehat{\ff}(\zeta)\right\|_\lt,
\]
where $\mathscr{S}'\left(\rr\right)$ is the space of tempered distributions.

Let $\x$ be the space defined by
\begin{equation}\label{x-space}
\x=\left\{f\in H^{1/2}(\rr);\, (\xi^{-1}\what{f}(\xi))^\vee\in L^2(\rr)\right\}
\end{equation}
 with the norm
 \[
{ \|f\|}_\x=\|f\|_{H^{1/2}(\rr)}+\|\partial_x^{-1}f\|_{L^2(\rr)}.
 \]
\section{Solitary Waves}

By a solitary wave solution of the RGBO equation, we mean a traveling-wave solution of equation \eqref{rgbo} of the form
$\ff(x-ct)$, where $\ff\in\x$ and $c\in \rr$ is the speed of wave propagation. Alternatively, it is a solution $\ff(x)$ in $\x$ of the stationary equation
\begin{equation}\label{rgbo-1}
\beta\h\ff_x-c\ff+f(\ff)=\ga\partial_x^{-2}\ff.
\end{equation}
We will prove existence of solitary waves in the space $\x$ by considering the following variational problem.
Define the functionals
\begin{equation}
I(u)=I(u;\beta,c,\ga)=\int_\rr\beta(\dx u)^2-cu^2+\ga(\dn u)^2\dd x
\end{equation}
and
\begin{equation}
K(u)=-(p+1)\int_\rr F(u)\dd x;
\end{equation}
and consider the following minimization problem
\begin{equation}\label{minimization-1}
M_\lam=\inf\{I(u);u\in\x, K(u)=\lam\},
\end{equation}
for some $\lam>0$.

First we observe that $M_\lam>0$ for any $\lam>0$. In fact, for $c\leq0$
\begin{equation}\label{coe-1}
\max\{\beta,c,\ga\}\|u\|_\x^2\geq I(u)\geq\beta\int_\rr(\dx u)^2\dd x+\ga\int_\rr(\dn u)^2\dd x;
\end{equation}
while for $c\in(0,c_\ast)$
\begin{equation}\label{coe-2}
\max\{\beta,c,\ga\}\|u\|_\x^2\geq I(u)\geq c_1\beta\int_\rr(\dx u)^2\dd x+c_2\ga \int_\rr(\dn u)^2\dd x,
\end{equation}
where $c_1=1-\sqrt{\frac{1}{2}+\frac{2c^3}{27\ga\beta^2}}$ and $c_2=\frac{27\beta^2-4c^3}{27\gamma\beta^2+4c^3}$.
On the other hand, for $u\in\x$, we have
\begin{equation}\label{emb}
\|u\|_{L^{p+1}(\rr)}^{p+1}\leq C{\|u\|}_{H^{1/2}(\rr)}^{(3p+1)/3}{\|\dn u\|}_{L^2(\rr)}^{2/3}.
\end{equation}
Indeed, by using the Sobolev embedding and an interpolation, we find
\begin{equation}\label{emb-1}
\|u\|_{L^{p+1}(\rr)}^{p+1}\leq C{\|u\|}_{H^{\frac{p-1}{2(p+1)}}(\rr)}^{p+1}\leq C{\|u\|}_{H^{1/2}(\rr)}^{(3p-1)/3}{\|u\|}_{H^{-1/4}(\rr)}^{4/3}.
\end{equation}
Then the Cauchy-Schwarz inequality implies that
\begin{equation}\label{emb-2}
\|u\|_{H^{-1/4}(\rr)}\leq C \|u\|^{1/2}_{H^{1/2}(\rr)}\|\dn u\|_{L^2(\rr)}^{1/2}.
\end{equation}
Now inequality \eqref{emb} is obtained from \eqref{emb-1} and \eqref{emb-2}.

We also note that $\|u\|_\x^2$ is equivalent to $I(u)$. Indeed, \eqref{coe-1}, \eqref{coe-2} and the inequality
\begin{equation}
\|u\|_\lt^2\leq A\|\dx u\|_\lt^2+B\|\dn u\|_\lt^2,\q\mbox{where}\q A>0,\;B>\frac{4}{27A^2},
\end{equation}
lead to the coercivity condition $I(u)\sim\|u\|_\x^2$.

Thus, it can be deduced from \eqref{emb} that
\[\begin{split}
\lam=K(u)\leq C\int_\rr|u|^{p+1}\dd x
&\leq
C{\|u\|}_{H^{1/2}(\rr)}^{(3p+1)/3}{\|\dn u\|}_{L^2(\rr)}^{2/3}\\
&\leq C\left(\|u\|_{L^2(\rr)}^2+\|\dx u\|_{L^2(\rr)}^2+\|\dn u\|_{L^2(\rr)}^2\right)^{(p+1)/2}.
\end{split}\]
Therefore we have $\lam^{2/(p+1)}\leq CI(u)$, where $C=C(\beta,c,\ga)>0$. This implies that $$M_\lam\geq(\lam/C)^{(p+1)/2}>0.$$

Then if $\psi\in\x$ achieves the minimum of problem \eqref{minimization-1}, for some $\lam>0$, then there exists a Lagrange multiplier $\mu\in\rr$ such that
\[
\beta\h\psi_x-c\psi-\ga\partial_x^{-2}\psi=-\mu f(\psi).
\]
Hence $\ff=\mu^{1/(p-1)}\psi$ satisfies \eqref{rgbo-1}.
We
denote the set of
such solutions by $G(\beta,c,\ga)$. By the homogeneity of $I$ and $K$, $u\in G(\beta,c,\ga)$ also achieve the minimum
\[
m=m(\beta,c,\ga)=\inf\left\{\frac{I(u)}{(K(u))^{\frac{2}{p+1}}};u\in\x, K(u)>0\right\}
\]
and it follows that $M_\lam=m\lam^{\frac{2}{p+1}}$. We note that if we multiply \eqref{rgbo-1} by $\ff$ and integrate, we find that $I(\ff;\beta,c,\ga)=K(\ff)$. Thus we may characterize the set $G(\beta,c,\ga)$ as
\[
G(\beta,c,\ga)=\left\{\ff\in\x; K(\ff)=I(\ff;\beta,c,\ga)=(m(\beta,c,\ga))^{\frac{p+1}{p-1}}\right\}.
\]
We now seek to prove that this set is nonempty.

We say that a sequence $\psi_n$ is a minimizing sequence if for some $\lam> 0$,
$\lim_{n\to\infty}K(\psi_n)=\lam$ and $\lim_{n\to\infty}I(\psi_n)=M_\lam$.
\begin{theorem}\label{exist}
Let $p\geq2$, $\beta> 0$, $\gamma>0$ and $c < c_\ast=3(\beta^2\ga/4)^{1/3}$. Let $\{\psi_n\}_n$ be a minimizing
sequence for some $\lam > 0$. Then there exist a subsequence (renamed $\psi_n$) and scalars $y_n\in\rr$ and $\psi\in\x$ such that $\psi_n(\cdot + y_n)\to\psi$  in $\x$. The function $\psi$ achieves the minimum $I(\psi) = M_\lam$ subject to the constraint $K(\psi) =\lam$.
\end{theorem}

\begin{figure}
\scalebox{0.3}{\includegraphics{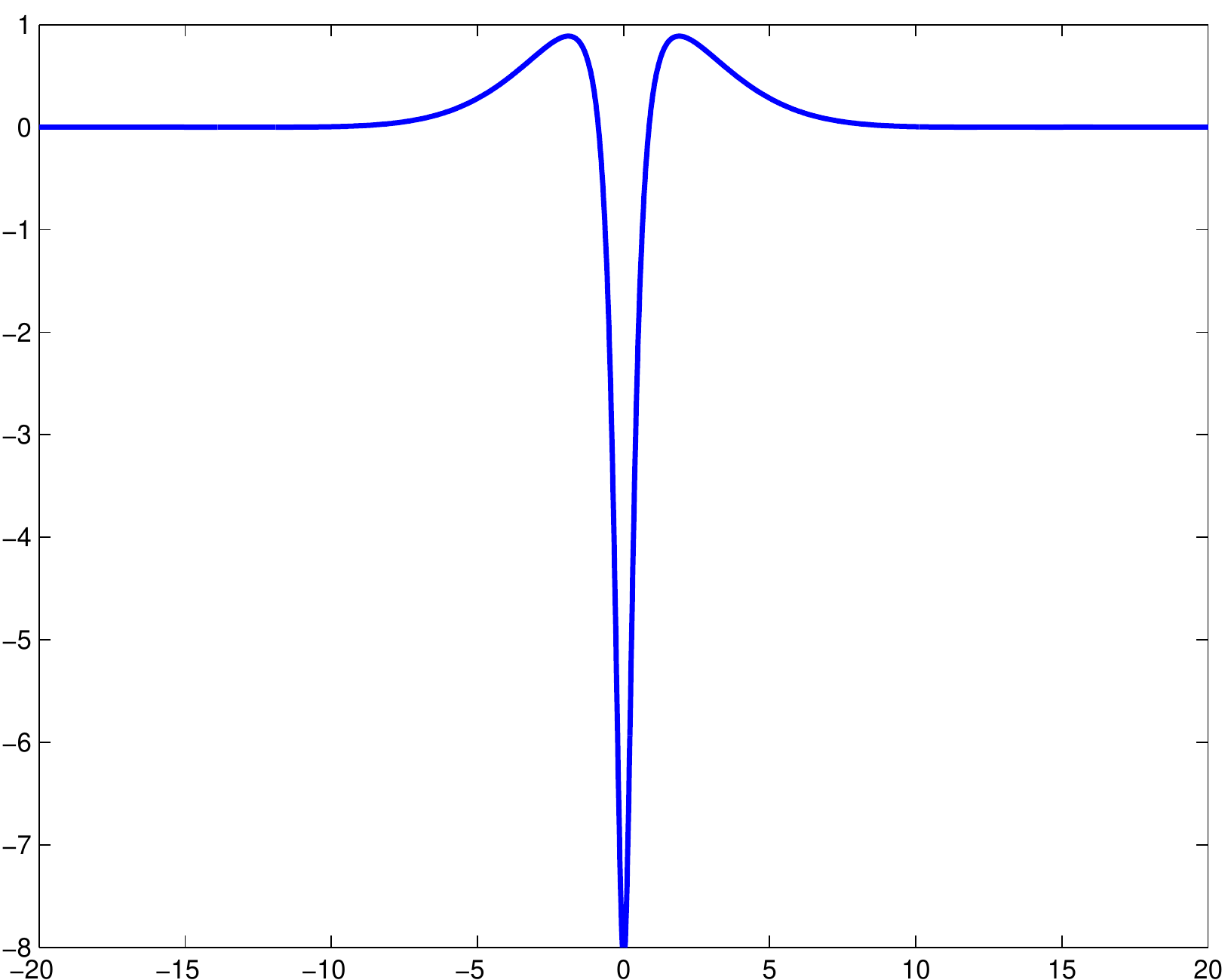}\quad
\includegraphics{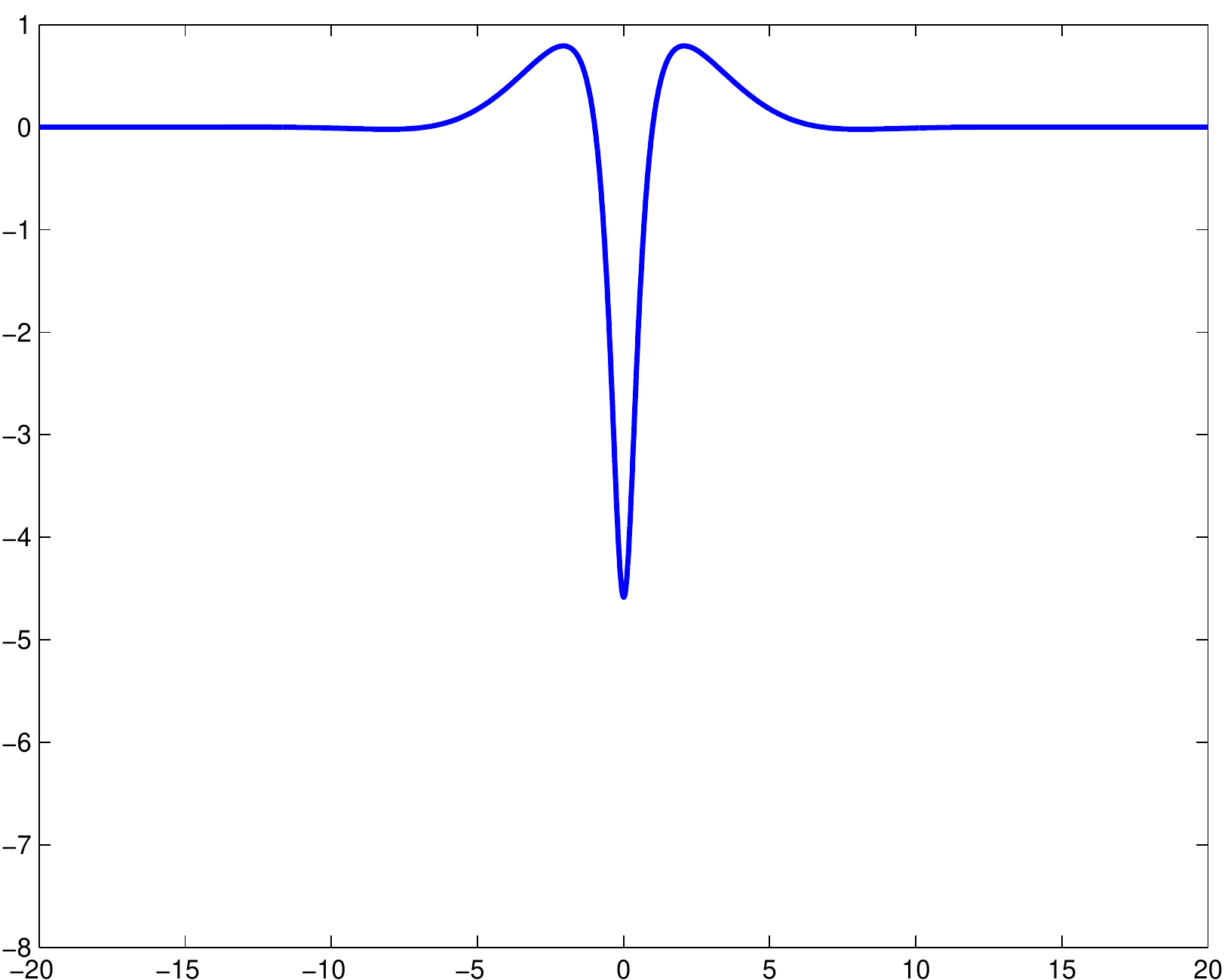}\quad
\includegraphics{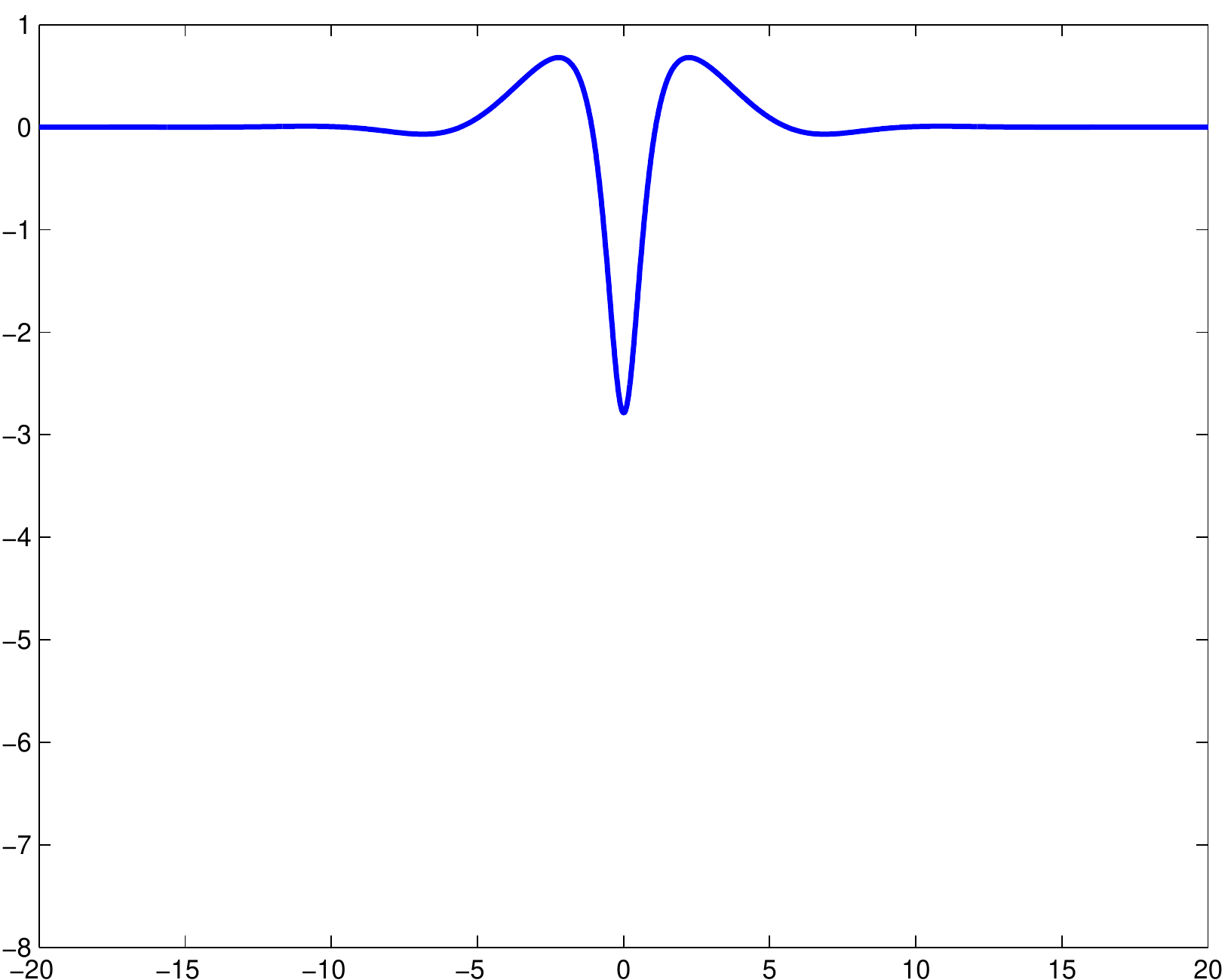}}
\scalebox{0.3}{
\includegraphics{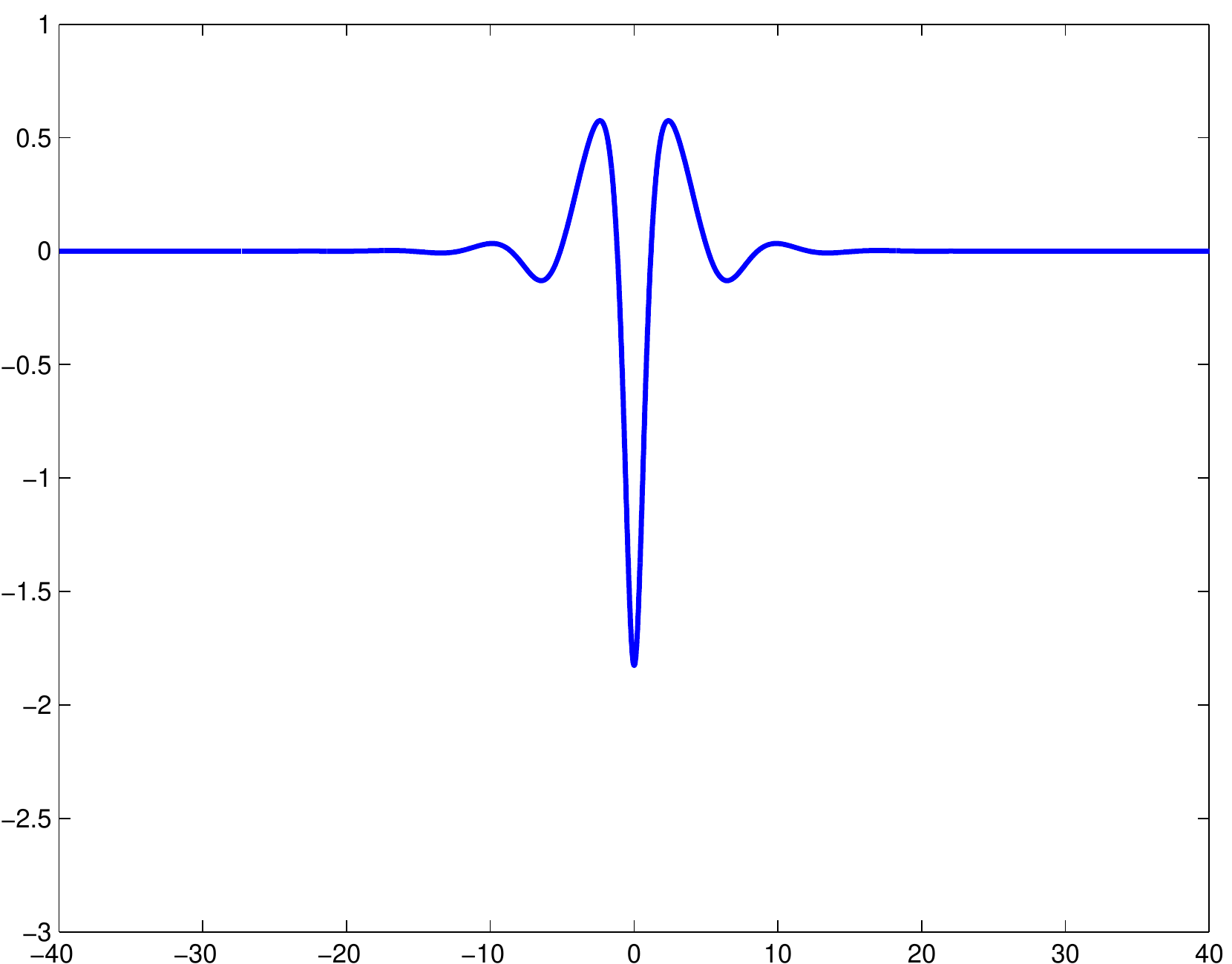}\quad
\includegraphics{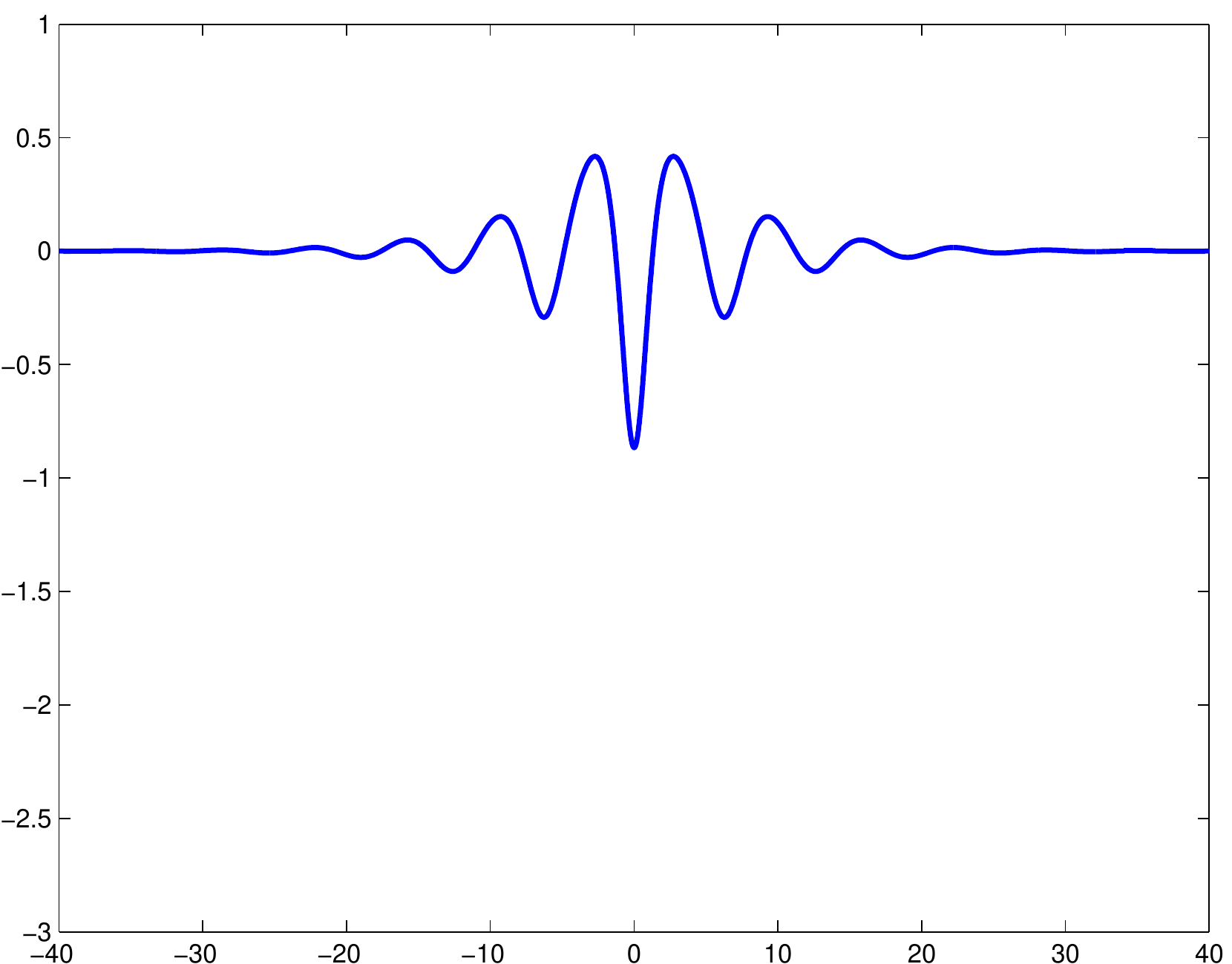}\quad
\includegraphics{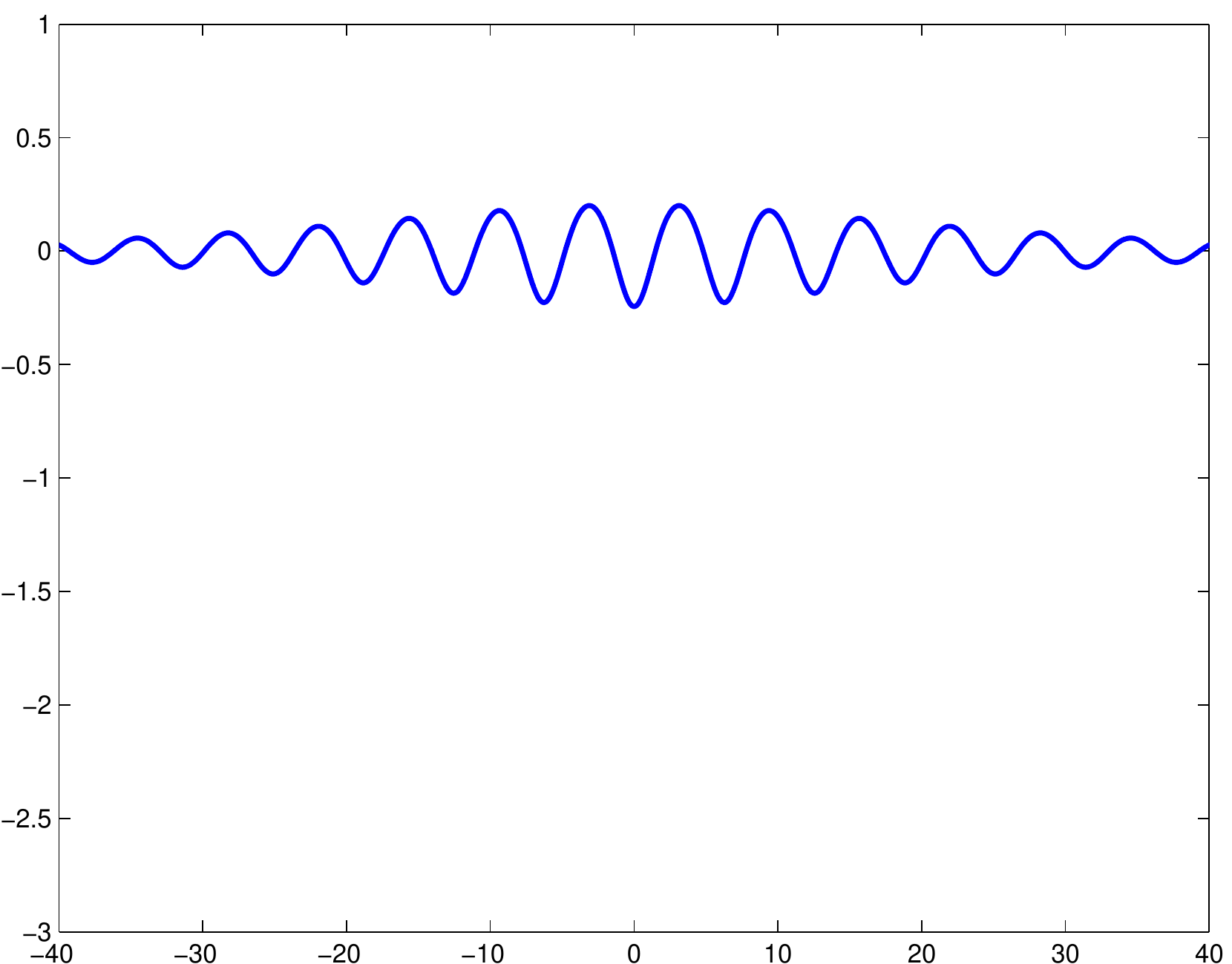}
}\caption{Solitary waves of the rotation generalized Benjamin-Ono equation with $f(u)=u^2$, $\beta=2$, $\gamma=1$ and $c=-1$, $1$, $2$, $2.5$, $2.9$ and $2.99$. Note that $c_\ast=3$ when $\beta=2$ and $\gamma=1$.}\label{F:c_limit}
\end{figure}

\proof
The result is an application of the concentration compactness lemma of Lions \cite{lions}, similar to \cite{albert,abs,liuvar}. We give the sketch of the proof here.

Let $\{\psi_n\}$ be a minimizing sequence, then we deduce from the coercivity of $I$ that the sequence $\{\psi_n\}$ is bounded in $\x$, so if we define
\[
\rho_n=|\dx\psi_n|^2+|\dn\psi_n|^2,
\]
then after extracting a subsequence, we may assume $\lim_{n\to\infty}\|\rho_n\|_{L^1(\rr)}=L>0$. We may assume further after normalizing that $\|\rho_n\|_{L^1(\rr)} = L$ for all $n$. By the concentration
compactness lemma, a further subsequence $\rho_n$ satisfies one of vanishing, dichotomy or compactness conditions. We can easily see that $M_\lam=\lam^{2/(p+1)}M_1$, so that the strict subadditivity condition
\[
M_\alpha+M_{\lam-\alpha}>M_\lam
\]
holds for any $\alpha\in(0,\lam)$. In the same manner as in \cite{leva,levaliu,liuvar}, it follows from the coercivity of $I$, inequality \eqref{emb}, and the subadditivity condition that both vanishing and dichotomy may be ruled out, and therefore the sequence $\rho_n$ is compact. Now set $\ff_n(x)=\psi(x+y_n)$. Since $\ff_n$ is bounded in $\x$, a subsequence $\ff_n$ converges weakly to some $\psi\in\x$, and
by the weak lower semicontinuity of $I$ over $\x$, we have $I(\psi)\leq\lim_{n\to\infty}I(\ff_n)=M_\lam$. Moreover, weak convergence in $\x$, compactness of $\rho_n$, and inequality \eqref{emb} imply strong convergence of $\ff_n$ to $\psi$ in $L^{p+1}(\rr)$. Therefore $K(\psi) = \lim K(\ff_n)=\lam$,
so $I(\psi)\geq M_\lam$. Together with the inequality above, this implies $I(\psi) = M_\lam$, so $\psi$ is a
minimizer of $I$ subject to the constraint $K(\cdot) = \lam$. Finally, since $I$ is equivalent to
the norm on $\x$, $\ff_n\rightharpoonup\psi$, and $I(\ff_n)\to I(\psi)$, it follows that $\ff_n$ converges strongly to
$\psi$ in $\x$.
\fim
\begin{theorem}\label{regularity}
Let $\beta$, $\ga$  and $c$   be as in Theorem \ref{exist}, and $f\in C^k$, for some nonnegative integer $k$. Then any weak solution $\ff$ of \eqref{rgbo-1} is a $H^{k+1}(\rr)$-function. Moreover $\dn\ff\in H^{k+2}(\rr)$ and $\partial_x^{-2}\ff\in H^{k+3}(\rr)$.
\end{theorem}
\proof
First we write \eqref{rgbo-1} in the form of a convolution equation
\begin{equation}\label{conv}
\ff=h\ast g(\ff),
\end{equation}
where $g(\ff)=-f(\ff)$ and
\begin{equation}\label{kernel}
\what{h}(\xi)=\frac{\xi^2}{\beta\xi^2|\xi|-c\xi^2+\ga}.
\end{equation}
Since $\frac{\xi^3}{\beta\xi^2|\xi|-c\xi^2+\ga}$ and $\frac{1}{\beta\xi^2|\xi|-c\xi^2+\ga}$ are bounded for any $\xi\in\rr$, then \eqref{conv} and the Sobolev embedding implies that $\partial_x^{-2}\ff,\ff_x\in L^2(\rr)$. Therefore we find that $\ff\in H^1(\rr)$; so that $\ff\in L^\infty(\rr)$. Hence $(g(\ff))_x\in L^2(\rr)$. Since
\[
(\beta\h\ff_x-c\ff-g(\ff))_x=\gamma\dn\ff
\]
it follows that $\h\ff_{xx}\in\lt$ and consequently $\ff_{xx}\in\lt$ and $\ff\in H^2(\rr)$. Repeating this process yields the property $\ff\in H^{k+1}(\rr)$.

A similar argument shows that $\dn\ff\in H^{k+2}(\rr)$ and $\partial_x^{-2}\ff\in H^{k+3}(\rr)$.
\fim

\begin{theorem}\label{analy}
Let $p>1$ be an integer, $\beta$, $\ga$  and $c$  be as in Theorem \ref{exist} and $\ff$  be a nontrivial solution of \eqref{rgbo-1}. There exist  $\kappa>0$ and an holomorphic function $\psi$ of  variable $z$, defined in the domain
\[
\mathcal{H}_\kappa=\left\{z\in\mathbb{C};\;|\mathrm{Im}(z)|<\kappa\right\},
\]
such that $\psi(x)=\ff(x)$ for all $x\in\rr$.
\end{theorem}
\proof
By the Cauchy-Schwarz inequality, we have that $\what{\ff}\in L^1\left(\rr\right)$. The equation \eqref{rgbo-1} implies that
\begin{equation}
\hspace{2mm}\left|\what{\ff}\right|(\xi)\leq
\overbrace{|\what{\ff}|\ast\cdots\ast |\what{\ff}|}^{p}(\xi),
\end{equation}
\begin{equation}
|\xi|\left|\what{\ff}\right|(\xi)\leq \underbrace{|\what{\ff}|\ast \cdots\ast|\what{\ff}|}_{p}(\xi),
\end{equation}
We denote $\mathscr{T}_1(|\what{\ff}|)=|\what{\ff}|$ and for
$m\geq1$,
$\mathscr{T}_{m+1}(|\what{\ff}|)=\mathscr{T}_m(|\what{\ff}|)\ast|\what{\ff}|$.
It can be easily seen by induction that for all $m\in\mathbb{N}$,
\begin{equation}
|\xi|^m|\what{\ff}|(\xi)\leq (m-1)!\;p^{m-1}\mathscr{T}_{(m+1)(p-1)+1}(|\what{\ff}|)(\xi).
\end{equation}
Therefore we have
\[\begin{split}
|\xi|^m|\what{\ff}|(\xi)
&\leq (m-1)!\;p^{m-1}{\left\|\mathscr{T}_{(m+1)(p-1)+1}(|\what{\ff}|)\right\|}_{L^\infty(\rr)}\\
&\leq (m-1)!\;p^{m-1}{\left\|\mathscr{T}_{(m+1)(p-1)}(|\what{\ff}|)\right\|}_{L^2(\rr)}{\left\|\what{\ff}\right\|}_{L^2(\rr)}\\
&\leq (m-1)!\;p^{m-1}{\left\|\what{\ff}\right\|}_{L^1(\rr)}^{(m+1)(p-1)-1}{\|\what{\ff}\|}_{L^2(\rr)}^2.
\end{split}\]
Let $$a_m=\frac{\;p^{m-1}
{\|\what{\ff}\|}_{L^1(\rr)}^{(m+1)(p-1)-1}{\|\what{\ff}\|}_{L^2(\rr)}^2}{m},$$
then it is clear that
\[
\frac{a_{m+1}}{a_m}\longrightarrow(p+1)\|\what{\ff}\|_{L^1(\rr)}^{p-1},
\]
as $m\to+\infty$. Therefore the series $\sum_{m=0}^\infty \zeta^m
r^m|\what{\ff}|(\xi)/m!$ converges uniformly in $L^\infty(\rr)$, if
$0<\zeta<\kappa=\frac{1}{(p)}\|\what{\ff}\|_{L^1(\rr)}^{-p+1}$. Hence
$\ee^{\zeta r}\what{\ff}(\xi,\eta)\in L^\infty(\rr)$, for
$\zeta<\kappa$.\\\indent We define the function
\[
\psi(z)=\int_{\rr}\ee^{\ii\xi z}\what{\ff}(\xi)\;\dd\xi.
\]
By the Paley-Wiener Theorem, $\psi$ is well defined and
analytic in $\mathcal{H}_\kappa$; and by Plancherel's Theorem, we
have $\psi(x)=\ff(x)$ for all $x\in\rr$. This proves the theorem.
\fim
\begin{theorem}\label{decay}
Let $\beta$, $\ga$   and $c$   be as in Theorem \ref{exist}. Then any solution $\ff$ of \eqref{rgbo-1} satisfies $|x|^\ell\ff^{(k)}(x)\in L^q(\rr)$, for $1\leq q\leq\infty$,  $k\in\{-2,-1,0\}$ and $\ell\in[0,4+k]$.  Furthermore,
\begin{equation}\label{decay-form}
|x|^\ell\ff^{(n)}(x)\in L^q(\rr),\q\mbox{for}\q 1\leq q\leq\infty,\; n\in\N,\; \ell\in[0,5].
\end{equation}
\end{theorem}
\proof
First a straightforward calculation reveals that $\what{h}\in C^\infty(\rr\setminus\{0\})\cap C^4(\rr)$. Moreover $\partial_\xi^j\what{h}\in L^q(\rr)$, for $q\in[1,\infty]$ and $1\leq j\leq4$. This implies that $\what{h}\in H^4(\rr)$. Hence by \cite[Corollary 3.1.3]{bonali}, we see that $\ff\in L^1(\rr)\cap\lt$ and $|x|^\ell\ff(x)\in\lt\cap L^\infty(\rr)$, for $\ell\in[0,4]$. Now the elementary inequality
\[
|x|^\ell|\ff|\leq ||x|^\ell h\ast g(\ff)|+| h\ast |x|^\ell g(\ff)|
\]
and the Young inequality imply that $|x|^\ell\ff(x)\in L^1(\rr)$, for $\ell\in[0,4]$.

Analogously, by using \eqref{rgbo-1}, one can show for $k=-2,-1$ that  $|x|^\ell\ff^{(k)}(x)\in L^q(\rr)$, for any $1\leq q\leq\infty$ and $\ell\in[0,4+k]$.

To prove \eqref{decay-form}, first we note that the fact $\ff\in L^\infty(\rr)$, the inequality
\[
|x|^\ell|\ff'|\leq ||x|^\ell h\ast (g(\ff))_x|+| h\ast |x|^\ell (g(\ff))_x|
\]
and the Young inequality implies that $|x|^\ell\ff'(x)\in L^q(\rr)$, for any $1\leq q\leq\infty$ and $\ell\in[0,4]$. On the other hand, a straightforward computation shows that $h'\in L^1(\rr)$ and $|x|^\ell h'\in L^q(\rr)$ for any $\ell\in[1,5]$ and $1\leq q\leq\infty$. Therefore combining the inequality
\[
|x|^5|\ff'|\leq ||x|^5 h'\ast g(\ff)|+| h'\ast |x|^5 g(\ff)|,
\]
the identity $|x|^5|\ff|^{p}=|x||\ff|(|x|^{4/(p-1)}|\ff|)^p$ and the Young inequality yields that $|x|^5\ff'(x)\in L^q(\rr)$, for any $1\leq q\leq\infty$. Finally a bootstrapping argument proves \eqref{decay-form}.
\fim

\begin{proposition}\label{6-decay-proposition}
Let $\beta>0$, $\ga>0$ and $c<c_\ast$ be as in Theorem \eqref{exist}, then there exists $C\in\rr$, $C\neq0$, such that any solution of \eqref{rgbo-1} satisfies
\begin{equation}\label{6-decay}
\lim_{|x|\to+\infty}|x|^6\ff(x)=C.
\end{equation}
\end{proposition}
\proof
The kernel $h$ in \eqref{kernel} can be written in $h(x)=-\frac{\dd^2}{\dd x^2}\kb(x)$, where
\[
\what{\kb}(\xi)=\frac{1}{\beta|\xi|^3-c\xi^2+\ga}.
\]
Since $\kb$ is an even function, hence
\begin{equation}\label{kernel-1}
\kb(x)=\int_\rr\frac{\ee^{\ii x\xi}}{\beta|\xi|^3-c\xi^2+\ga}\dd\xi=\int_0^{+\infty}\frac{\cos(|x|\xi)}{\beta\xi^3-c\xi^2+\ga}\dd\xi.
\end{equation}
Then by using the residue theorem, there holds that
\begin{equation}
\begin{split}
\kb(x)&=
\int_0^{+\infty}\frac{-\beta y^3\ee^{-y|x|}}{(\ga+cy^2)^2+\beta^2 y^6}\dd y
+2\pi
{\rm Re}\left(\frac{\ii\ee^{b\ii|x|-a|x|}}{3\beta(b^2-a^2)-2cb)+2a\ii(3b\beta-2c)}\right)\\
&=
\int_0^{+\infty}\frac{-\beta y^3\ee^{-y|x|}}{(\ga+cy^2)^2+\beta^2 y^6}\dd y\\
&\q+2\pi
\frac{\ee^{-a|x|}\left[2a(3b\beta-c)\cos(bx)+\sin(b|x|)(2cb-3\beta(b^2-a^2))\right]}
{\left(3\beta(b^2-a^2)-2cb\right)^2+4a^2\left(c-3b\beta\right)^2},
\end{split}
\end{equation}
where $b+\ii a$ is the complex root of $\beta|\xi|^3-c\xi^2+\ga$, with $a,b>0$.
Therefore $\kb\in C^\infty(\rr\setminus\{0\})$. It is therefore  concluded for $c<c_\ast$ that
\begin{equation}\label{h-kernel}
\begin{split}
h(x)&=\int_0^{+\infty}\frac{\beta y^5\ee^{-y|x|}}{(\ga+cy^2)^2+\beta^2 y^6}\dd y\\
&\q-2\pi a(3b\beta-c)\frac{\ee^{-a|x|}\left(2ab\sin(b|x|)+(a^2-b^2)\cos(bx)\right)}
{\left(3\beta(b^2-a^2)-2cb\right)^2+4a^2\left(c-3b\beta\right)^2}\\
&\q-2\pi(2cb-3\beta(b^2-a^2))\frac{\ee^{-a|x|}\left((a^2-b^2)\sin(b|x|)-2ab\cos(bx)\right)}
{\left(3\beta(b^2-a^2)-2cb\right)^2+4a^2\left(c-3b\beta\right)^2}.
\end{split}\end{equation}
Now the change of variable $\eta=xy$ in the first term of the right-hand side of \eqref{h-kernel} reveals that
\[
\lim_{|x|\to+\infty}|x|^6h(x)=\frac{\beta}{\ga^2}.
\]
Applying Theorem 3.1.5 in \cite{bonali}, it transpires that there exists $C\neq0$ such that \eqref{6-decay} holds.
\fim
\begin{remark}
By Theorem \ref{decay} and Proposition \ref{6-decay-proposition}, one can see that the solitary waves of \eqref{rgbo-1} decay faster than the solitary waves of \eqref{gbo} (see \eqref{explicit}).
\end{remark}
\begin{remark}
By \eqref{h-kernel}, it seems that the solitary wave $\ff$ of \eqref{rgbo-1} does not decay exponentially.
\end{remark}
\section{Nonexistence}

In this section we present conditions on the parameters $\beta$, $c$, $\gamma$ and the nonlinearity $f(u)$ that guarantee equation \eqref{rgbo} has no solitary wave solutions in the space $\x$. These conditions follow from the following functional identities.
\begin{lemma} Suppose $\ff\in\x$ is a solution of equation \eqref{rgbo-1}. Then
\begin{equation}\label{E:functional_identities}
\begin{array}{lll}
  \int_\rr\beta(\dx\ff)^2\dd x-c\int_\rr \ff^2\dd x+\ga\int_\rr(\dn \ff)^2\dd x =-(p+1)\int_\rr F(\ff)\dd x\\
  \\
  -c\int_\rr \ff^2\dd x+3\ga\int_\rr(\dn \ff)^2\dd x =-2\int_\rr F(\ff)\dd x\\
\end{array}
\end{equation}
\end{lemma}

\proof These relations follow by multiplying equation \eqref{rgbo-1} by $\ff$ and $x\ff_x$, respectively and integrating over $\rr$. To see that the $\beta$ term vanishes in the second relation, first observe that since $\ff_x$ has zero mass, it follows that $\h(x\ff_x)=x\h(\ff_x)$. Then, using the anti-commutative property of $\h$ we have
    \[
    \int_\rr\h\ff_x\cdot x\ff_x\dd x
    =-\int_\rr \ff_x\cdot\h(x\ff_x)\dd x
    =-\int_\rr \ff_x\cdot x\h\ff_x\dd x.
    \]
This completes the proof.
\fim

\begin{theorem} Equation \eqref{rgbo-1} has no solution in $\x$ provided any of the following conditions hold.
\begin{enumerate}[(i)]
  \item $\beta<0$, $\gamma>0$ and $c^3<\frac{27(3p+1)\gamma\beta^2}{(p-1)^3}$.
  \item $\beta>0$, $\gamma<0$ and $c^3>\frac{27(3p+1)\gamma\beta^2}{(p-1)^3}$.
  \item $f(u)=|u|^{p-1}u$, $\beta>0$ and $\gamma<0$.
  \item $f(u)=-|u|^{p-1}u$, $\beta<0$ and $\gamma>0$.
\end{enumerate}
\end{theorem}

\proof Eliminating the terms on the right hand sides of \eqref{E:functional_identities}, we find that
    \begin{equation}\label{E:homogeneous_identity}
    -2\beta\int_\rr(\dx\ff)^2\dd x-(p-1)c\int_\rr \ff^2\dd x+(3p+1)\ga\int_\rr(\dn \ff)^2\dd x =0.
    \end{equation}
Now suppose $\beta<0$ and $\gamma>0$. Then since the expression
    \[
    (3p+1)\gamma|\xi|^{-2}-2\beta|\xi|
    \]
has minimum value $-3\beta(3p+1)^{1/3}(-\gamma/\beta)^{1/3}>0$ it follows that
    \[
    -2\beta\int_\rr(\dx\ff)^2\dd x+(3p+1)\ga\int_\rr(\dn \ff)^2\dd x
    \geq-3\beta(3p+1)^{1/3}(-\gamma/\beta)^{1/3}\int_\rr \ff^2\dd x,
    \]
so if $c$ satisfies the inequality in (i), the left hand side of \eqref{E:homogeneous_identity} will be negative, a contradiction. This proves statement (i). Statement (ii) follows similarly.

Next, subtracting the two relations in \eqref{E:homogeneous_identity}, we have
    \[
    -\beta\int_\rr(\dx\ff)^2\dd x+2\ga\int_\rr(\dn \ff)^2\dd x =(p-1)\int_\rr F(u)\dd x.
    \]
The right and left hand sides of this equation have opposite signs when either condition (iii) or condition (iv) holds.
\fim

\section{Ground States and Variational Characterizations}
A ground state of \eqref{rgbo-1} is a solitary wave of \eqref{rgbo} which minimizes the action $S(u) = E(u) - c Q (u)$
among all nonzero solutions of \eqref{rgbo-1}, where $E(u)$ and $Q(u)$ are defined in \eqref{cons-1} and $\eqref{cons-2}$,
respectively.
Recall that a solitary wave of \eqref{rgbo} corresponds to a
critical point of $S(u)$, that is, $S'(u) = 0$. Thus, the set of ground states may be characterized as
\begin{equation}
\g = \{\ff\in\x;\; S'(\ff) = 0, S(\ff) \leq S(\psi)\,\mbox{ for all}\; \psi\in\x\;\mbox{satisfying}\; S'(\psi) = 0\}.
\end{equation}
The theorem below finds a ground state of \eqref{rgbo-1} as a minimizer for $S(\ff)$ under a new constraint. Our result is related to that in \cite{liu-0}.
\begin{theorem}\label{ground}
If $\beta$, $c$ and $\ga$ are as in Theorem \ref{exist}, then $\g$ is nonempty and $\ff\in\g$ if and only if $S(\ff)$ solves the minimization problem
\begin{equation}\label{minimization-2}
J = \inf\{S(u);\; \psi\in\x, \psi\neq0,\; P(\psi) = 0\},
\end{equation}
where
\[
P(\psi)=\int_\rr(-c\psi^2+\beta(\dx\psi)^2+\ga(\dn\psi)^2+(p+1)F(\psi))\dd x.
\]
\end{theorem}
\proof
First, we prove that there is a nontrivial minimizer for \eqref{minimization-2} which is a solution of \eqref{rgbo-1}.

By \eqref{emb}, one can easily observe that there exist $\varepsilon_1$, $\varepsilon_2>0$ such that for every nontrivial function $\ff\in\mathscr{N}$, we have $\|\ff\|_\x\geq\varepsilon_1$ and $S(\ff)\geq\varepsilon_2$, where
$\mathscr{N}=\{\psi\in\x;\;u\neq0,\;P(\psi)=0\}$.

Now, let $\{\ff_n\}\subset\mathscr{N}$ be a minimizing sequence of \eqref{minimization-2}. Then  $\|\ff_n\|_\x\geq\varepsilon_1$ and
\[
S(\ff_n)=\frac{p-1}{2(p+1)}I(\ff_n)\cong\frac{p-1}{2(p+1)}\|\ff_n\|_\x^2,
\]
so that $\{\ff_n\}_n$ is bounded in $\x$. To show that there is a convergent subsequence, with a limit $\ff\in\x$, similar to \cite{albert,abs}, we use again the concentration-compactness lemma \cite{lions}, applied to the sequence
\[
\rho_n=|\dx\ff_n|^2+|\dn\ff_n|^2.
\]
First similar to Theorem \ref{exist}, the evanescence case is excluded. To rule out the dichotomy case, one shows that
\[
J<J_\sigma:=\inf\left\{S(\psi)-\frac{1}{2}P(\psi);\;P(\psi)=\sigma\right\},
\]
for all $\sigma<0$. Now if the dichotomy would occur, i.e. $\ff_n$ splits into a sum $\ff_n^1+\ff_n^2$ and the distance of the supports of these functions tends to $+\infty$, then one shows that $P\left(\ff_n^1\right)\to\sigma$, $P\left(\ff_n^2\right)\to-\sigma$, $\sigma\in\rr$ and $J\geq J_\sigma+J_{-\sigma}>J$ which is a contradiction. Therefore the sequence $\ff_n$ concentrates and the limit $\ff$ satisfies $P(\ff)\leq0$. The case $P(\ff)<0$ can be excluded by the same reason as above, and we see that $\ff\in\mathscr{X}$ is a minimizer for \eqref{minimization-2}.

Now since $\ff$ is a minimizer for \eqref{minimization-2}, there exists a Lagrange multiplier $\theta$  such that $S'(\ff)=\theta P'(\ff)$. Since $\left\langle S'(\ff), \ff\right\rangle=0$ and
\[
\left\langle P'(\ff),\ff\right\rangle=2I(\ff)-(p+1)K(\ff)=(1-p)I(\ff)<0,
\]
we see that $\theta=0$, i.e. $\ff$ is a solution of \eqref{rgbo-1}.

Our next task is to show that $\ff\in\g$. But since $P(u)=\langle S'(u),u\rangle_{\lt}$ for any $u\in\x$, it follows that $P(v)=0$ for any solitary wave $v\in\x$ of \eqref{rgbo-1}. Hence $S(\ff)=J$ asserts that $S(\ff)\leq S(v)$.

Finally we show that a ground state of \eqref{rgbo-1} achieves the minimum $J$ in \eqref{minimization-2}. Let
$u\in\x$ satisfy $u\neq0$, $S'(u)=0$ and $S(u)\leq S(v)$ for any $v\in\x$ satisfying $S'(v) = 0$. Since
$S'(v) = 0$ implies $P(v) = \langle S'(v), v\rangle_\lt = 0$, it follows that $S(u)\leq S(v)$ for any $v\in\x$ with
$P(v) = 0$. That is, $v$ is a minimizer for $J$. This completes the proof.
\fim
The following proposition proves that minima for $M_\lam$ in \eqref{minimization-1} are exactly the ground states
of \eqref{rgbo-1}.
\begin{proposition}
There is a positive real number $\lam^\ast$ such that the following statements are equivalent:
\begin{enumerate}[(i)]
\item $K(\ff)=\lam^\ast$ and $\ff$ is a minimizer of $M_{\lam^\ast}$ in \eqref{minimization-1};
\item $\ff$ is a ground state;
\item $P(\ff) = 0$ and $K(\ff) = \inf\{K(u); u \in \x, u \neq 0,\; P(u) = 0\}$;
\item $P(\ff) = 0 = \inf\{P(u);\;u\in\x\; u\neq0,\;K(u) = K(\ff)\}$.
\end{enumerate}
\end{proposition}
\proof
(ii)$\Rightarrow$(i). Let $\ff$ be a ground state of \eqref{rgbo-1}. Since $P(\ff) = I(\ff) - K (\ff) = 0$ and $S(\ff) = \frac{1}{2}I(\ff) - \frac{1}{p+1}K (\ff)$, it follows that $\ff$ minimizes $I$ among solutions of \eqref{rgbo-1}. Set $\lam^\ast=K(\ff)=I(\ff)$.

Let $v$ be a minimizer for $M_{\lam^\ast}$. That is, $K(v) =\lam^\ast$ and $I(v) =M_{\lam^\ast}$ minimizes $I(u)$
among $K(u) =\lam^\ast$. In particular $M_{\lam^\ast}=I(v)\leq I(\ff)=\lam^\ast$. From variational considerations, $v$
satisfies
\[
-cv+\beta\h v_x-\ga\partial_x^{-2}v=-\theta f(v),
\]
for some $\theta\in\rr$. Multiplication of the above by $v$ and integration by parts yields
$I(v)=\theta K(v)$. Since $I(v)=M_{\lam^\ast}$, this implies $\theta\leq1$. On the other hand, since $w = \theta^{\frac1{p-1}}v$ is a solution of \eqref{rgbo-1}, one obtains $\theta^{\frac2{p-1}}I(v)=I(w)\geq I(\ff)$. Since $I(v)=\theta\lam^\ast$ and $I(\ff)=K(\ff) =\lam^\ast$, this implies $\theta\geq1$. Therefore, $\theta= 1$ and $I(\ff) =M_{\lam^\ast}$.

(i)$\Rightarrow$(iii). Suppose $K(\ff)=\lam^\ast$ and $\ff\in\x$ is a minimizer for $M_{\lam^\ast}$. Note that $P(\ff)=0$ and $I(\ff)=K(\ff)=\lam^\ast$. Let $u\in\x$ be such that $u\neq0$ and $P(u)=0$. Then $K(u)\neq>0$ so we may define $b=(K(\ff)/K(u))^{1/(p+1)}$. We show that $b\leq1$.

Straightforward calculations yield that $P(bu)=b^2(1-b^{p-1})I(u)$. Since $K(bu)=b^{p+1}K(u)=K(\ff)=\lam^\ast$, it follows that $I(\ff)\leq I(bu)$, and consequently $0=P(\ff)=I(\ff)-K(\ff)\leq I(bu)-K(bu)=b^2(1-b^{p-1})I(u)$. Hence the assertion follows.

 (ii)$\Leftrightarrow$(iii) is  a direct consequence of Theorem \ref{ground}.

(iii)$\Rightarrow$(iv). Let $u\in\x$, $u\neq 0$ with $K(u) = K(\ff)$,
where $\ff\in\x$ satisfies (iv). We prove that $P(u)\geq0$. Suppose on the contrary that
$P(u) < 0$. Note that $P(\tau u) > 0$ for $\tau\in(0, 1)$ sufficiently small. Correspondingly,
$K(u) > 0$ must hold and $P(\tau_0u) = 0$ for some $\tau_0\in(0, 1)$. This however contradicts (iii)
since $K(\tau_0 u)<K(u)=K(\ff)$. Therefore, $P(u)\geq 0$. The assertion then follows since $P(\ff) = 0$.

(iv)$\Rightarrow$(iii).
Let $u \in\x$, $u \neq0$ with $P(u) = 0$.
We show that $K(u)\geq K(\ff)$, where $\ff\in\x$ satisfies (iv). Assume the
opposite inequality. Similarly as in the previous argument, a scaling consideration shows that
$K (\tau_0u) = \tau^{p+1}_0K(u) = K(\ff)$, for some $\tau_0>1$. This contradicts (iii) since $P(\tau_0u) < 0 = P(\ff)$. This
completes the proof.
\fim

\section{Weak Rotation Limit}
In this section, we show that the solitary waves of the RGBO equation \eqref{rgbo} converge to those of the generalized Benjamin-ono equation \eqref{gbo}.
We remark that such a relationship is somewhat surprising since the solitary waves of \eqref{rgbo} have zero mass, as can be seen by integrating \eqref{rgbo} with respect to $x$, while it is well-known (see \cite{albert,albert-1,albert-2,albona,abr} and references therein) that the  solitary waves of \eqref{gbo} are strictly negative functions and do not have zero mass.

In order to precisely state the convergence result, it is worth noting that for each $c<0$ and $\beta>0$ the GBO equation \eqref{gbo} possesses a nontrivial solitary wave $\ff$  and it satisfies
\begin{equation}\label{solitary-gbo}
-c\ff+\beta\h\ff'+f(\ff)=0.
\end{equation}
The uniqueness of solitary waves of the GBO equation for $p>1$ is unknown; however Amick and Toland \cite{amicktoland} showed that the solitary wave solutions of the classical Benjamin-Ono ($p=1$) are unique (up to translation). The explicit solution was found by Benjamin \cite{benjamin}:
\begin{equation}\label{explicit}
\ff(\xi)=\frac{4c\beta^2}{\beta^2+c^2\xi^2}.
\end{equation}
One can see that contrary to the unique solitary wave of the KdV equation, the solitary wave of the Benjamin-Ono equation does not decay exponentially.

\begin{theorem}\label{weak}
For $\beta> 0$ and $c < 0$ fixed, let a sequence $\ga_n\to0^+$ as $n\to\infty$, and let $\psi_n$ any element of $G(\beta,c,\ga_n)$. Then there exists a subsequence (still denoted as $\ga_n$), translations $y_n$ and a solitary wave $\psi\in\st$ of \eqref{solitary-gbo} so that $\psi_n(\cdot+y_n)\to\psi$ in $H^{1/2}(\rr)$, as $\ga_n\to0^+$. That is, the solitary waves of the GBO equation are the limits in $H^{1/2}(\rr)$ of solitary waves of the  RGBO equation.
\end{theorem}

\begin{figure}[ht]
\scalebox{0.3}{\includegraphics{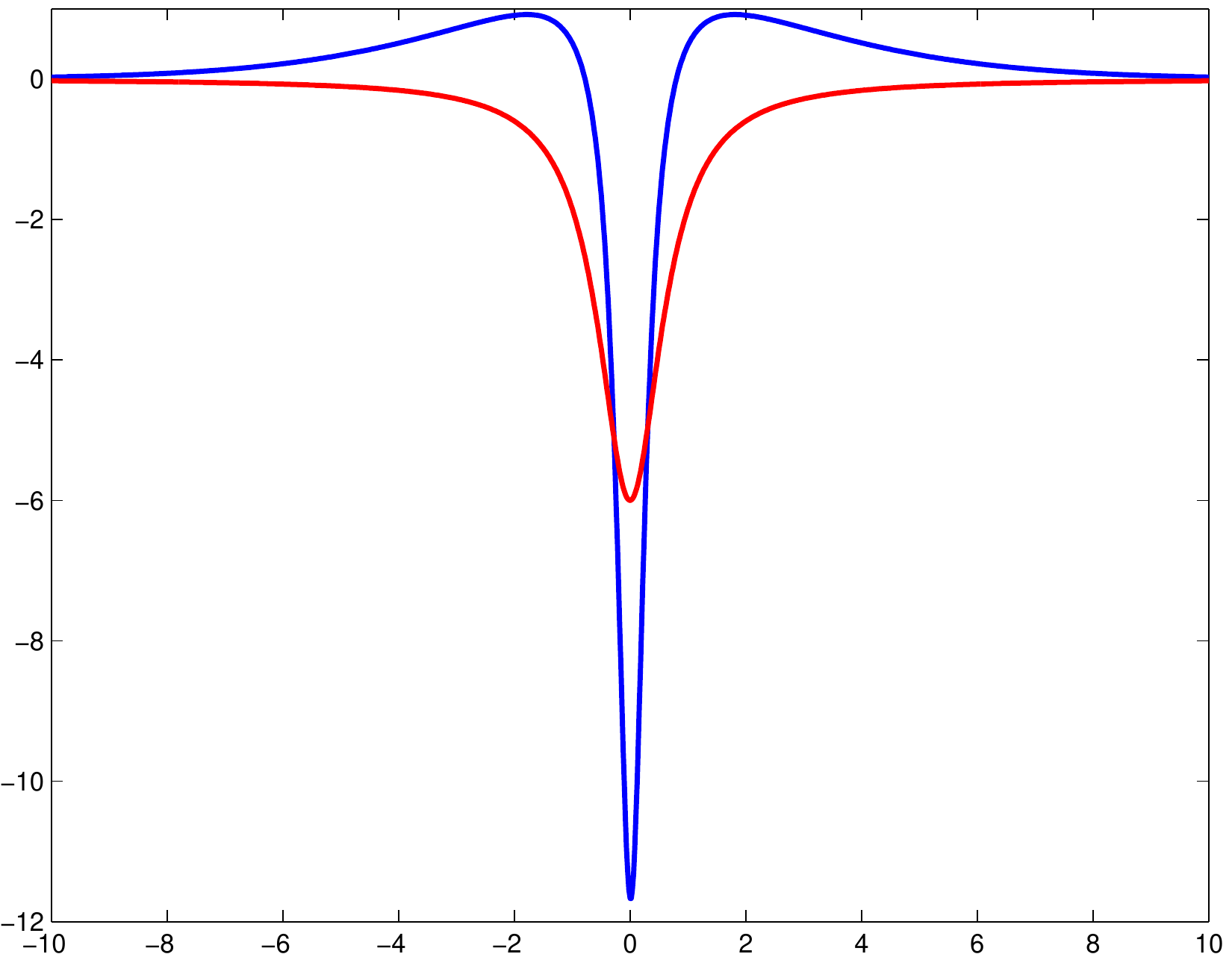}\quad
\includegraphics{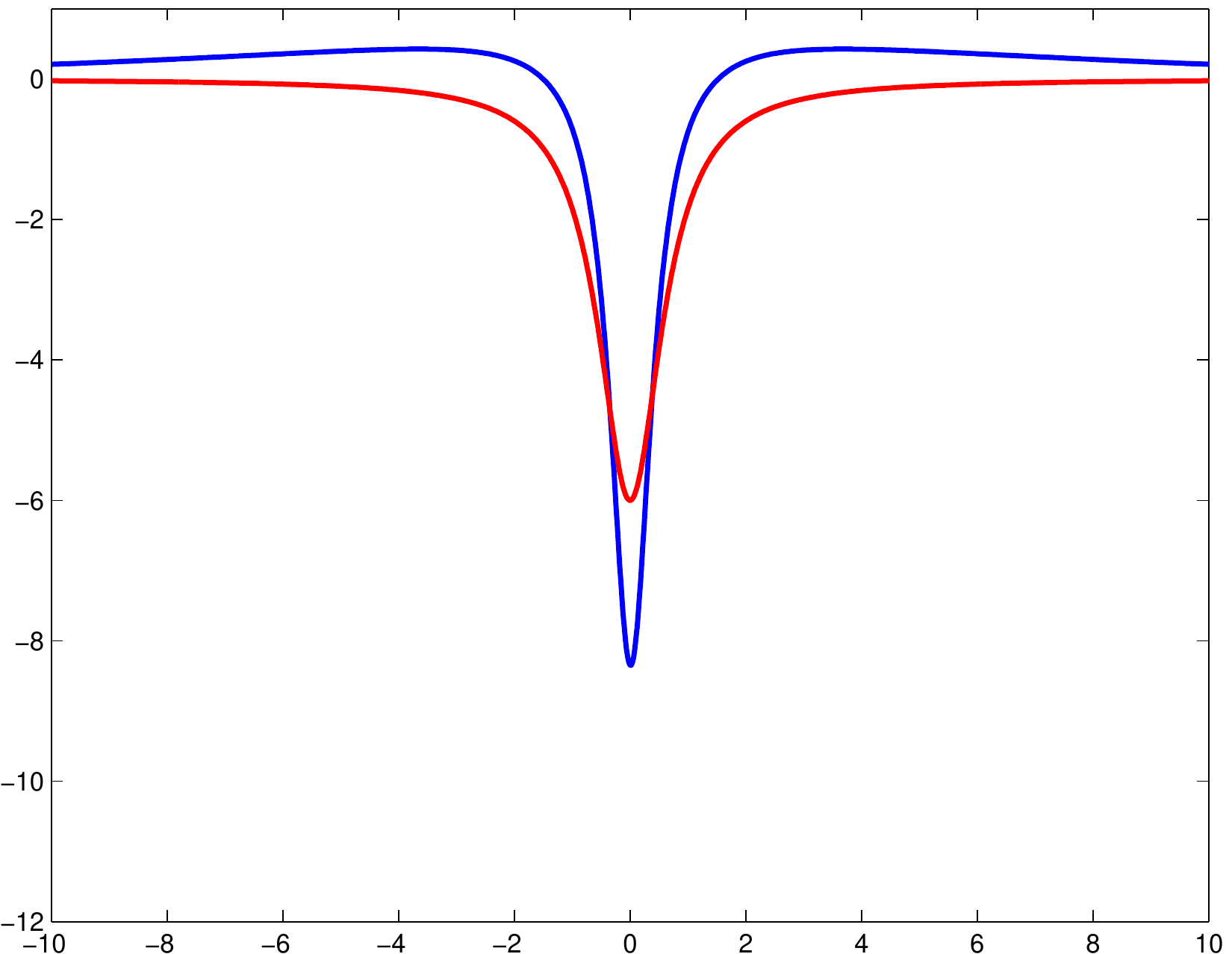}\quad
\includegraphics{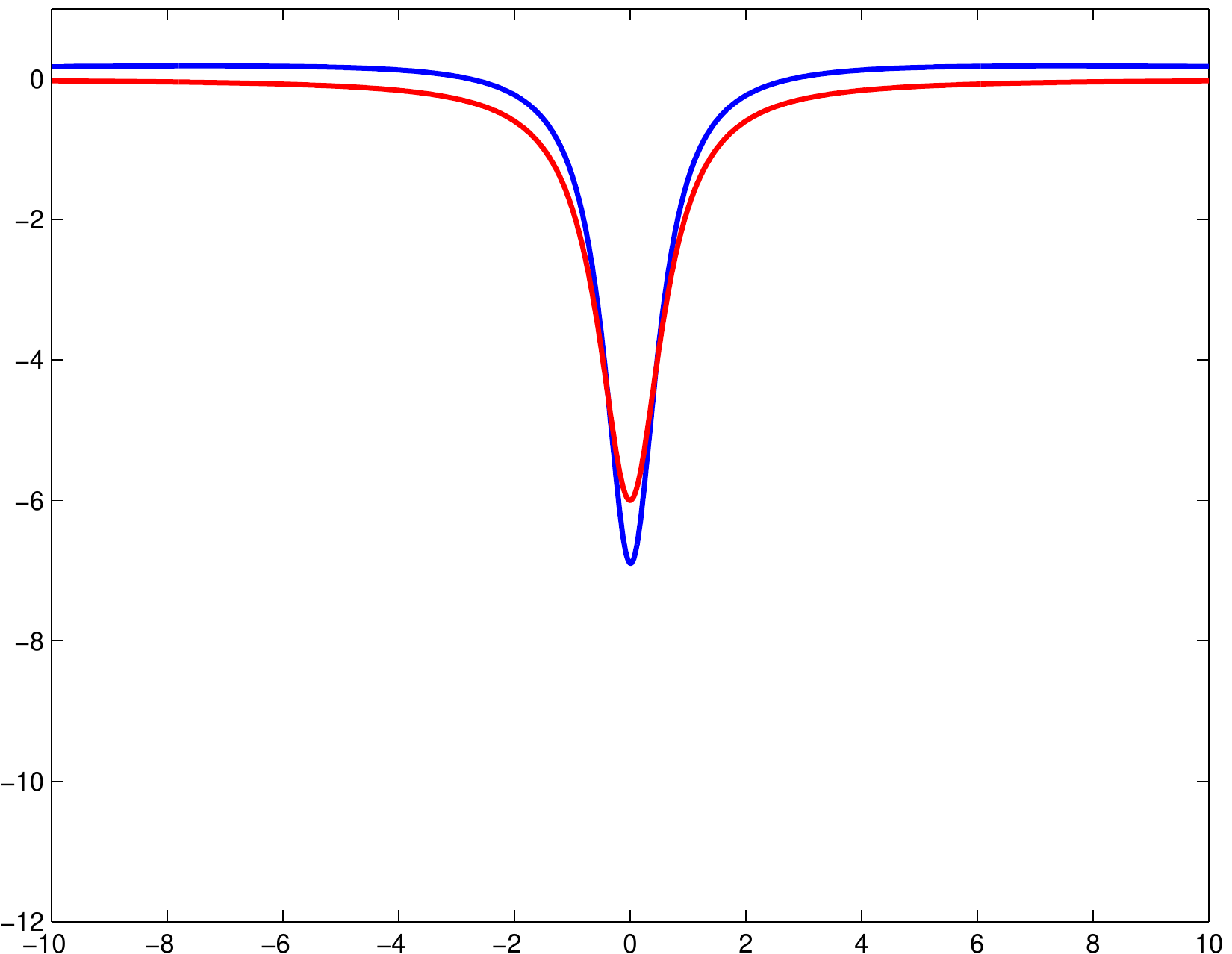}}
\scalebox{0.3}{\includegraphics{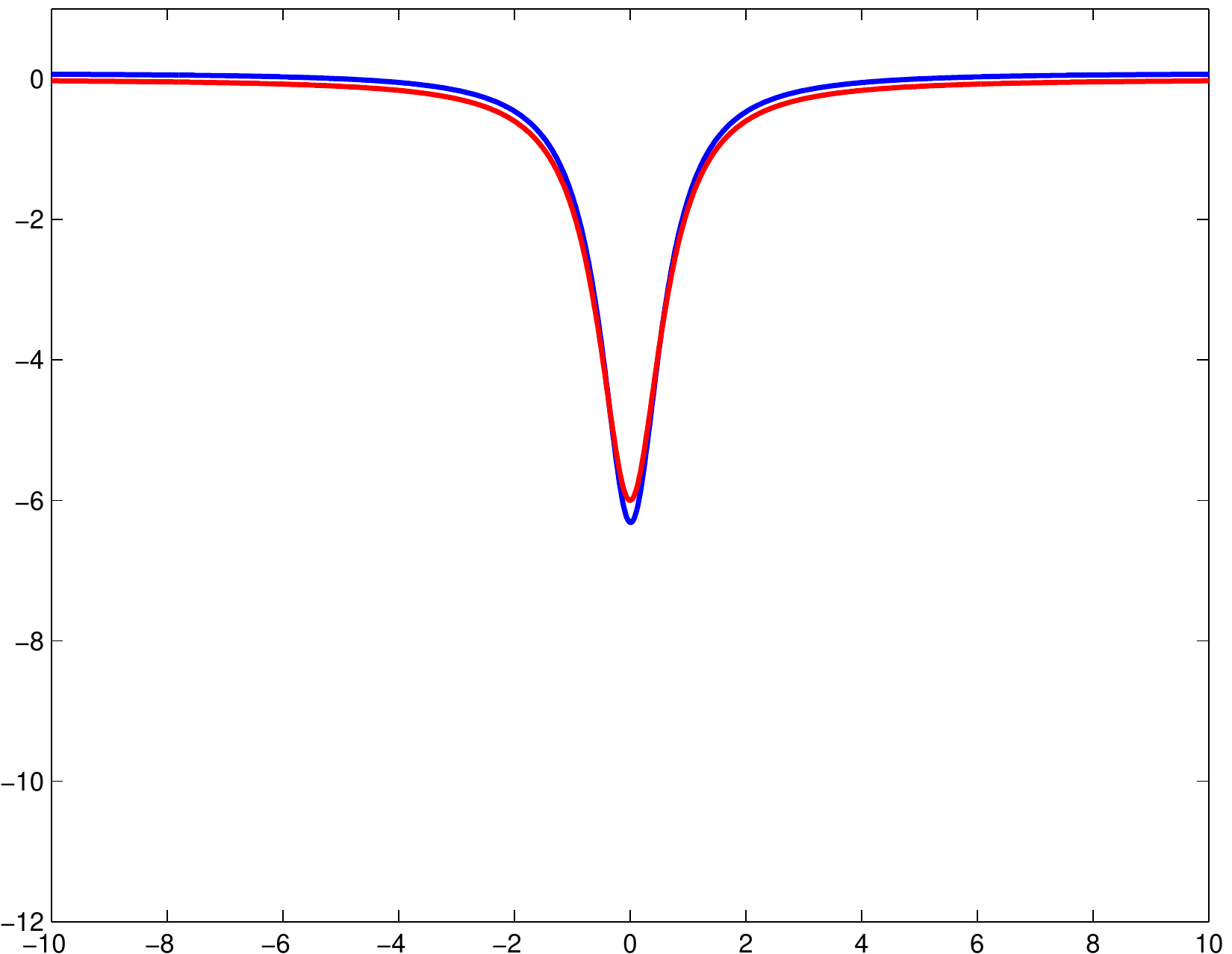}\quad
\includegraphics{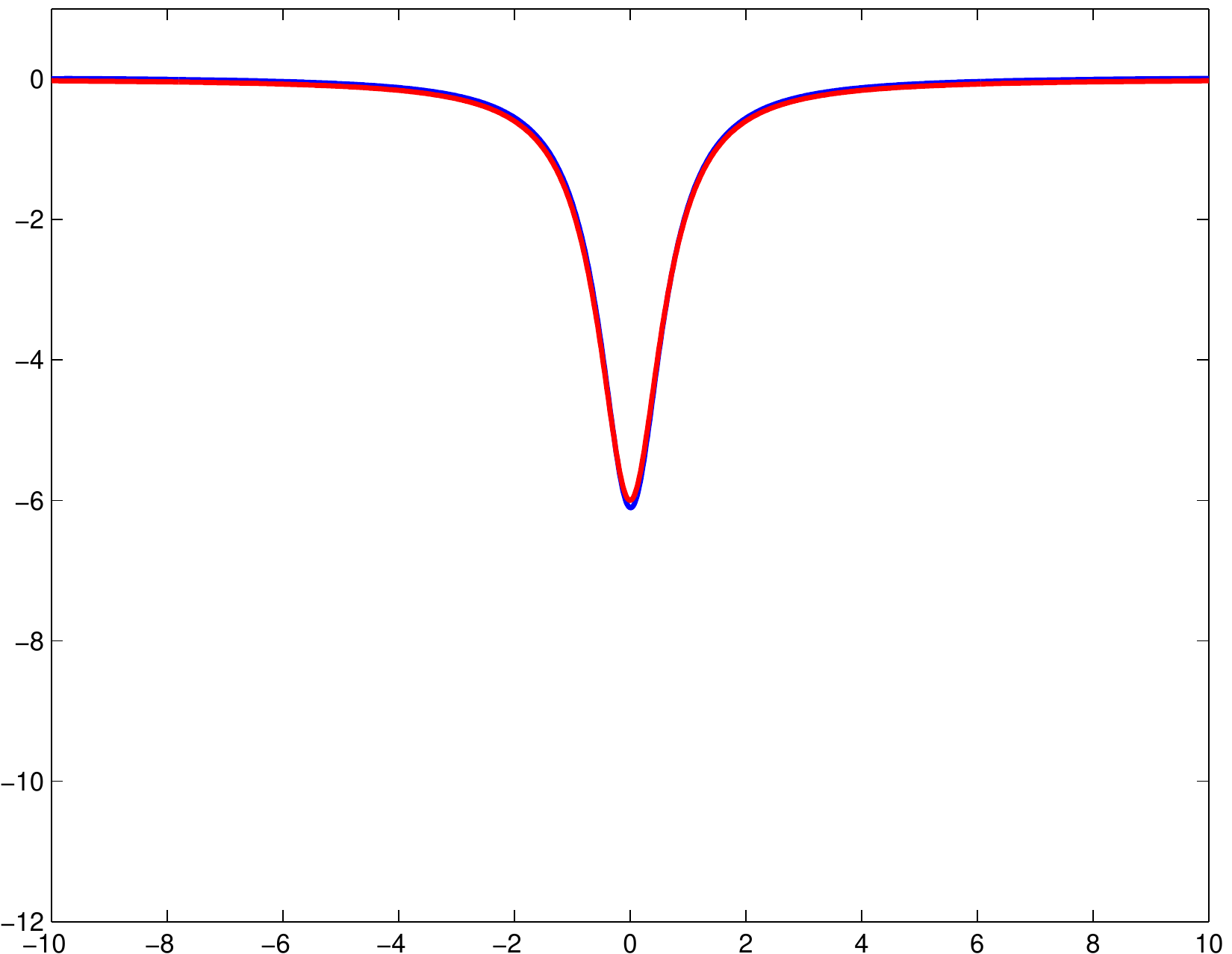}\quad
\includegraphics{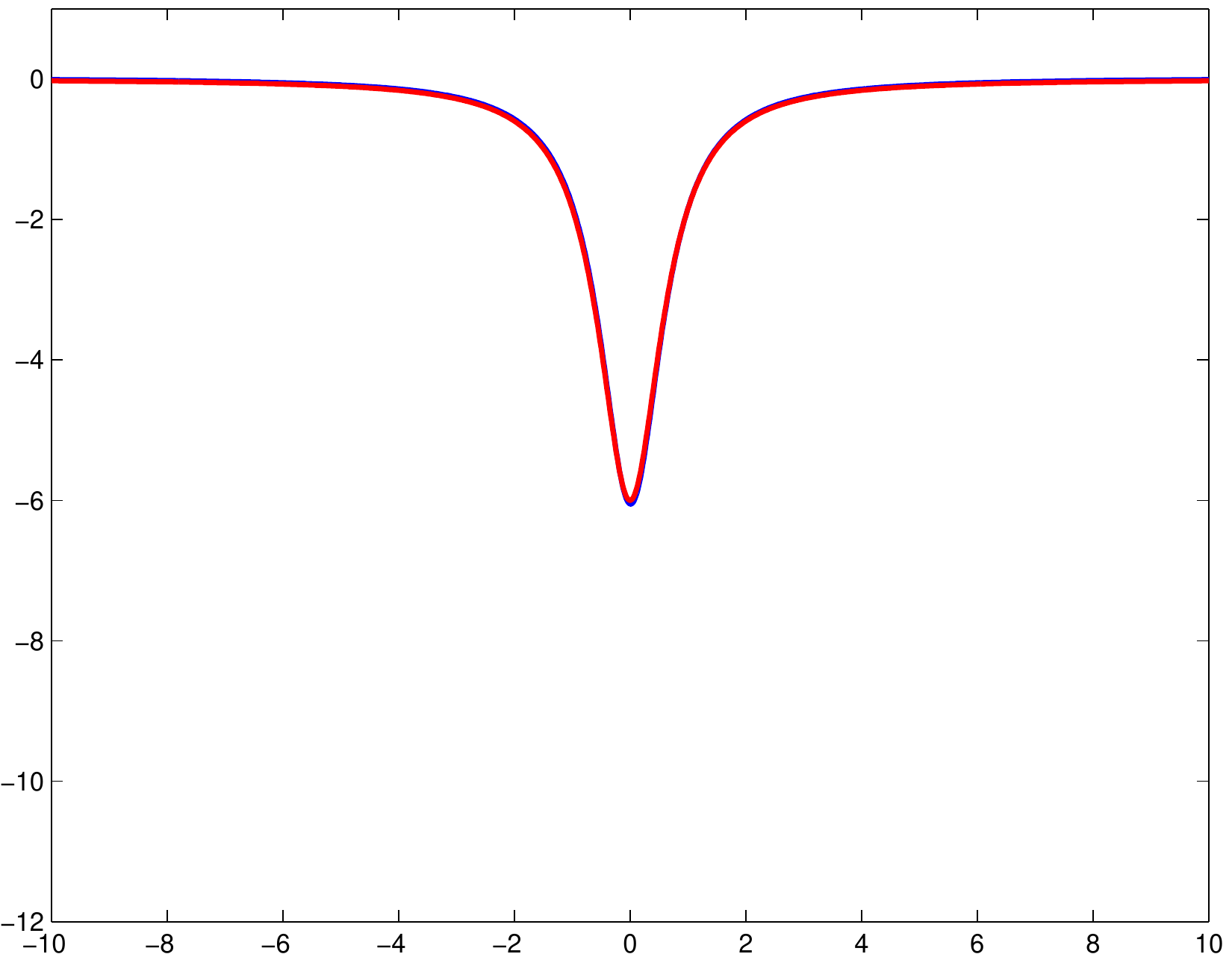}}
\caption{Solitary waves of the rotation generalized Benjamin-Ono equation with $f(u)=u^2$, $\beta=2$, $c=-3$ and $\gamma=1$, $0.1$, $0.01$, $0.001$, $0.0001$ and $0.00001$ are shown in blue. The exact solitary wave solution of the Benjamin-Ono given by \eqref{explicit} is shown in red.}\label{F:gamma_limit}
\end{figure}

To prove this result, we first note for $\beta>0$ and $c<0$ that solutions of \eqref{solitary-gbo} satisfies in a variational problem of the type of Theorem \ref{exist}. More precisely,  ground states of \eqref{solitary-gbo} achieve the minimum
\[
 m(\beta,c,0)=\inf\left\{\frac{I(u;\beta,c,0)}{K(u)^{\frac{2}{p+1}}}; u\in H^{1/2}(\rr),\; K(u)>0\right\},
\]
where $I(u;\beta,c,0)=\int_\rr(\beta(\dx u)^2-cu^2)\dd x$. Analogous to Theorem \ref{exist}, one can show that for a given sequence $\psi_n$ in $H^{1/2}(\rr)$ satisfying
\[
\lim_{n\to\infty}I(\psi_n;\beta,c,0)=\lim_{n\to\infty}K(\psi_n)=(m(\beta,c,0))^{\frac{p+1}{p-1}},
\]
there exists a subsequence, renamed $\psi_n$,  scalars $y_n\in\rr$ and $\ff_0\in H^{1/2}(\rr)$ such that $\psi_n(\cdot+y_n)\to\ff_0$ in $H^{1/2}(\rr)$.

The proof of Theorem \ref{weak} is approached via the following lemmas.

\begin{lemma}\label{monoton}
The function $m$ is continuous on the domain $\beta>0$, $\gamma>0$, $c<c_\ast$. Furthermore, $m$ is strictly increasing in $\ga$ and $\beta$ and strictly decreasing in $c$.
\end{lemma}

\proof
The proof is similar to \cite[Lemma 2.3]{levaliu}, \cite[Lemma 2.4]{levaliu-2} and  \cite[Lemma 3.3]{liuvar} by using the following inequality
\begin{equation}\label{best-estimate}
I(u;\beta,c,\ga)\geq (c_\ast-c)\int_\rr u^2\dd x,
\end{equation}
where $c_\ast$ is defined in Theorem \ref{exist}.
\fim

\begin{lemma}\label{density}
The space $\x$ is dense in $\st$.
\end{lemma}

\proof
For any $u\in\st$ and $\delta > 0$, let us define $u_\delta$ as
$\what{u_\delta}(\xi)=\what{u}(\xi)\chi_{|\xi|>\delta}(\xi)$. By Parseval's identity follow that
\[
\|\dn u_\delta\|_\lt^2=\|\xi^{-1}\what{u_\delta}\|_\lt^2=\int_{|\xi|>\delta}\xi^{-2}|\what{u}(\xi)|^2\dd\xi<
\delta^{-2}\|u\|_\lt^2<+\infty.
\]
Since $\|u_\delta\|_\lt\leq\|u\|_\lt<+\infty$ and since $\|\dx u_\delta\|_\lt\leq\|\dx u\|_\lt<+\infty$, it follows that $u_\delta\in\x$. In view of the definition of $u_\delta$ and $u\in\st$ then the inequality
\[
\|u_\delta-u\|_\st^2=\int_{|\xi|<\delta}(1+|\xi|)|\what{u}(\xi)|^2\dd\xi\leq\|u\|_\st^2<+\infty
\]
holds true. Hence from continuity we may choose $\delta>0$ sufficiently small so that
\[
\|u_\delta-u\|_\st^2=\int_{|\xi|<\delta}(1+|\xi|)|\what{u}(\xi)|^2\dd\xi<\epsilon,
\]
which completes the proof.
\fim

\noindent\textbf{Proof of Theorem \ref{weak}.}\quad
For $\beta > 0$ and $c < 0$ let $\{\psi_n\}$ be a sequence in $\x$ of the ground states of $\eqref{rgbo-1}$ with $\ga=\ga_n$, where $\ga_n\to0^+$ as $n\to\infty$ It is immediate that $I(\psi_n;\beta,c,\ga)=K(\psi_n)=m(\beta,c,\ga_n)^{\frac{p+1}{p-1}}$ holds for each $n$. Below we prove the continuity of $m(\beta,c,\ga)$ at $\ga=0$, that is, $\lim_{\ga\to0^+}m(\beta,c,\ga)=m(\beta,c,0)$. The assertion then follows from
\[
I(\psi_n;\beta,c,0)=I(\psi_n;\beta,c,\ga_n)-\ga_n{\|\dn\psi_n\|}_\lt^{2}
\leq I(\psi_n;\beta,c,\ga_n)=m(\beta,c,\ga_n)^{\frac{p+1}{p-1}}\longrightarrow m(\beta,c,0)^{\frac{p+1}{p-1}}
\]
and
\[
K(\psi_n)=m(\beta,c,\ga_n)^{\frac{p+1}{p-1}}\longrightarrow m(\beta,c,0)^{\frac{p+1}{p-1}}.
\]
We now claim that $\lim_{\ga\to0^+}m(\beta,c,\ga)=m(\beta,c,0)$. By the monotonicity of $m(\beta,c,\ga)$ in $\ga$ , it suffices to show that $m(\beta,c,\ga_n)\to m(\beta,c,0)$ for some sequence $\{\ga_n\}$ with $\ga_n\to0$ as $n\to\infty$. Let $\ff\in\st$ be a ground state of \eqref{solitary-gbo}. For each $n$ a positive integer it follows from
lemma \ref{density} that there is a function $\psi_n\in\x$ with $\|\psi_n-\ff\|_\st<1/n$. Let
\[
\ga_n=\min\left\{\frac{1}{n},\frac{1}{n}\|\dn\psi_n\|_\lt^{-2}\right\}.
\]
Then $\gamma_n\to0$ and
\[
m(\beta,c,\ga_n)\leq\frac{I(\psi_n;\beta,c,\gamma)}{K(\psi_n)^{\frac{2}{p+1}}}=
\frac{I(\psi_n;\beta,c,0)+\ga_n{\|\dn\psi_n\|}_\lt^{2}}{K(\psi_n)^{\frac{2}{p+1}}}
\leq\frac{I(\psi_n;\beta,c,0)+1/n}{K(\psi_n)^{\frac{2}{p+1}}}.
\]
Since both $I(\cdot;\beta,c,0)$ and $K$ are continuous on $\st$, it follows that
\[
\lim_{n\to\infty}m(\beta,c,\ga_n)\leq\frac{I(\ff;\beta,c,0)}{K(\ff)^{\frac{2}{p+1}}}=m(\beta,c,0).
\]
On the other hand, since $m(\beta,c,\ga)$ is strictly increasing in $\ga$, it follows that
\[
\lim_{n\to\infty}m(\beta,c,\ga_n)=m(\beta,c,0).
\]
This proves the claim. The proof is complete.
\fim

\section{Stability}
In this section we investigate the stability of the set $\g$ of ground state solitary waves. We first state precisely our definition of stability.
\begin{definition}
A set $\Omega\in\x$ is $\x$-stable with respect to \eqref{rgbo} if for any $\e>0$ there exists $\delta>0$ such that for any $u_0\in\x\cap X_s$, $s>3/2$, with
\[
\inf_{v\in \Omega}\|u_0-v\|_\x<\delta,
\]
then the solution $u(t)$ of \eqref{rgbo} with initial value $u(0) = u_0$ can be extended to a solution in the space $C([0,+\infty),\x\cap X_s)$ and satisfies
\[
\sup_{t\geq0}\inf_{v\in \Omega}\|u(t)-v\|_\x<\e.
\]
Otherwise we say that $\Omega$ is $\x$-unstable.
\end{definition}

Since ground state solitary waves minimize the action $S(u)=E(u)-cQ(u)$, it is natural to consider the function
\begin{equation}\label{d-function}
d(c)= d(\beta,c,\gamma)=E(\varphi)-cQ(\varphi),
\end{equation}
where $\ff$ is any element of $\g$. The fact that $d$ is well-defined follows from the relation
    \begin{equation}\label{E:dIKM_relation}
    d(\beta,c,\gamma)=\frac12I(\varphi)-\frac1{p+1}K(\varphi)=\frac{p-1}{2(p+1)}K(\varphi)
    =\frac{p-1}{2(p+1)}m(\beta,c,\gamma)^{\frac{p+1}{p-1}}.
    \end{equation}
Together with Lemma \ref{monoton}, this relation also implies that $d$ is continuous on the domain $\beta>0$, $\ga>0$, $c<c_\ast$, strictly increasing in $\ga$ and $\beta$ and strictly decreasing in $c$. It can also be shown as in \cite{levaliu} and \cite{levaliu-2} that $d$ has the following differentiability properties.

\begin{lemma}\label{L:d_derivatives} For each fixed $\beta>0$ and $\gamma>0$,
the partial derivative $d_c(\beta,c,\gamma)$
exists for all but countably many $c$. For fixed $c$ and $\gamma$, $d_\beta(\beta,c,\gamma)$ exists for all but countably many $\beta$ and for fixed $\beta$ and $c$, $d_\gamma(\beta,c,\gamma)$ exists for all but countably many $\gamma$. At points
of differentiability, we have
\begin{align*}
    d_\beta(\beta,c,\gamma)&=\frac12\int(D_x^{1/2}\varphi)^2\,\dd x\\
  d_c(\beta,c,\gamma)&=-\frac12\int \varphi^2\,\dd x=-Q(\varphi)\\
  d_\gamma(\beta,c,\gamma)&=\frac12\int (\partial_x^{-1}\varphi)\,\dd x
\end{align*}
\end{lemma}

For the remainder of this section we shall regard $\beta>0$ and $\gamma>0$ as fixed and denote $d(c)=d(\beta,c,\gamma)$, $d'(c)=d_c(\beta,c,\gamma)$ and $d''(c)=d_{cc}(\beta,c,\gamma)$.
The main stability result is that the stability of the set of ground states is determined by the sign of $d''(c)$.

\begin{theorem}\label{stability-theo}
Let $\beta>0$, $\gamma>0$, $c < c_\ast$ and $\ff\in\g$. If $d''(c)>0$, then the set of ground states $\g$ is $\x$-stable.
\end{theorem}

Define the $\epsilon$-neighborhood of the set of ground states defined by
\[
U_{\epsilon}=\{u\in\x;\;\inf_{\ff\in\g}\|u-\ff\|_\x<\epsilon\}.
\]
Since $d$ is strictly decreasing in $c$ and $K$ is continuous on $\x$, we may define
\[
c(u)=d^{-1}\left(\frac{p-1}{2(p+1)}K(u)\right)
\]
for $u\in U_\epsilon$ for sufficiently small $\epsilon>0$.
\begin{lemma}\label{taylor-exp}
If $d''(c)>0$, then there is some $\eps>0$ such that for any $u\in U_\eps$
and $\ff\in\g$, we have
\begin{equation}
E(u)-E(\ff)-c(u)(Q(u)-Q(\ff))\geq\frac{1}{4}d''(c)|c(u)-c|^2.
\end{equation}
\end{lemma}
\proof
Since $d'(c) = -Q (\ff)$, it follows from Taylor's theorem that
\[
d(c_1)=d(c)-Q(\ff)(c_1-c)+\frac{1}{2}(c_1-c)^2+o(|c_1-c|^2),
\]
for $c_1$ near $c$. Using the continuity of  $c(u)$ and choosing $\e>0$ sufficiently small we get that
\[
d(c(u))\geq d(c)-Q(\ff)(c(u)-c)+\frac{1}{4}d''(c)(c(u)-c)^2
=E(\ff)-c(u)Q(\ff)+\frac{1}{4}(c(u)-c)^2,
\]
for $u\in U_\eps$. Next, if $\ff_{c(u)}\in\mathscr{G}(\beta,c(u),\ga)$; then $K(\ff_{c(u)})=2(p+1)d(c(u))/(p-1)=K(u)$ and $\ff_{c(u)}$  minimizes $I(\cdot;\beta,c(u),\ga)$ subject to this constraint, so
\[
E(u)-c(u)Q(u)=\frac{1}{2}I(u;\beta,c(u),\ga)-\frac{1}{p+1}K(u)
\geq
\frac{1}{2}I(\ff_{c(u)};\beta,c(u),\ga)-\frac{1}{p+1}K(\ff_{c(u)})=d(c(u)).
\]
This concludes the proof of Lemma \ref{taylor-exp}.
\fim
\noindent\textbf{Proof of Theorem \ref{stability-theo}.}\quad
Assume that $\g$ is $\x$-unstable with regard to the flow of the RGBO equation. Then there exists a sequence of the initial data $u_k(0)$ such that
\[
\inf_{\ff\in\g}\|u_k(0)-\ff\|_\x<\frac{1}{k}.
\]
Let $u_k(t)$ be the solution of \eqref{rgbo} with initial data $u_k(0)$. We can also choose $\delta>0$ and a sequence of times $t_k$ such that
\begin{equation}\label{stab-1}
\inf_{\ff\in\g}\|u_k(t)-\ff\|_\x=\delta.
\end{equation}
Moreover we can find $\ff\in\g$ such that
\[
\lim_{k\to\infty}\|u_k(0)-\ff_k\|_\x=0.
\]
Since $E$ and $V$ are conserved by the flow of \eqref{rgbo},
\begin{equation}\label{stab-2}
\lim_{k\to\infty}E(u_k(t_k))-E(\ff_k)=\lim_{k\to\infty}E(u_k(0))-E(\ff_k)=0
\end{equation}
and
\begin{equation}\label{stab-3}
\lim_{k\to\infty}Q(u_k(t_k))-Q(\ff_k)=\lim_{k\to\infty}Q(u_k(0))-Q(\ff_k)=0.
\end{equation}
By using Lemma \ref{taylor-exp}, we have for $\delta$ sufficiently small that
\begin{equation}\label{stab-4}
E(u_k(t_k))-E(\ff_k)-c(u_k(t_k))\left(Q(u_k(t_k))-Q(\ff_k)\right)\geq\frac{1}{4}d''(c)|c(u_k(t_k))-c|^2.
\end{equation}
By \eqref{stab-1} there is some $\psi_k\in\g$ such that $\|u_k(t_k))\|_\x<2\delta$, and by using the fact $I(u)=I(u;\beta,c,\gamma)\geq C\|u\|_\x^2$, we obtain
\[
\|u_k(t_k))\|_\x\leq\|\psi_k\|_\x+2\delta\leq C^{-1}I(\psi_k;\beta,c,\gamma)+2\delta=\frac{2(p+1)}{C(p-1)}d(c)+2\delta<\infty.
\]
Thus since $K$ is Lipschitz continuous on $\x$ and $d^{-1}$ is continuous, it follows that
$c(u_k(t_k))$ is uniformly bounded in $k$. Thus by \eqref{stab-2}-\eqref{stab-4} it follows that $\lim_{k\to\infty}c(u_k(t_k))=c$; and therefore
\begin{equation}\label{stab-5}
\lim_{k\to\infty}K(u_k(t_k))=\lim_{k\to\infty}\frac{2(p+1)}{(p-1)}d(c(u_k(t_k)))=\frac{2(p+1)}{(p-1)}d(c).
\end{equation}
This implies that
\[\begin{split}
\frac{1}{2}I(u_k(t_k))&=E(u_k(t_k))-cQ(u_k(t_k))+\frac{1}{p+1}K(u_k(t_k))\\
&=
d(c)+E(u_k(t_k))-E(\ff_k)-c(Q(u_k(t_k))-Q(\ff_k))+\frac{1}{p+1}K(u_k(t_k)).
\end{split}\]
Hence it follows from \eqref{stab-2}, \eqref{stab-3} and \eqref{stab-5} that
$\lim_{k\to\infty}I(u_k(t_k))=2(p+1)d(c)/(p-1)$. Thus $u_k(t_k)$ is a minimizing sequence and therefore has a subsequence which converges in $\x$ to some $\ff\in\g$. This contradicts \eqref{stab-1}, so the proof of the theorem is complete.
\fim

\section{Instability}

In this section we present conditions that imply orbital instability of ground state solitary waves.

Given $\ff\in\g$ and $\epsilon>0$, we define
\[
\Omega_{\ff,\eps}=\left\{u\in\x;\;\inf_{v\in\mathcal{O}_\ff }\|v-u\|_\x<\eps\right\}.
\]
\begin{theorem}\label{instability-cond-theo}
Let $\beta>0$, $\ga>0$, $c<c_\ast$ and $\ff\in\g$. Suppose there exists $\psi\in\lt$ such that $\psi'\in X_s$, $s>3/2$, $\psi''\in\x$, and the following two conditions hold.
\begin{equation}\label{instability-cond}
\begin{split}
&\left<\psi',\ff\right\>=0,\\
&\left\<S''(\ff)\psi',\psi'\right\><0.
\end{split}
\end{equation}
Then $\mathcal{O}_\ff$ is $\x$-unstable.
\end{theorem}
\begin{lemma}\label{implicit}
Let $c<c_\ast$ and $\ff\in\g$ be fixed. There are an $\eps_0>0$ and a unique $C^2$ map $\al:\Omega_{\ff,\eps_0}\to\rr$ such that $\al(\ff)=0$, and for all $v\in \Omega_{\ff,\eps_0}$ and any $r\in\rr$,
\begin{enumerate}[(i)]
\item $\<\tau_{\al(v)}\ff',v\>=0$,
\item $\al(\tau_rv)=\al(v)+r$,
\item $\al'(v)=-\frac{1}{\left\<v,\ff''(\cdot+\al(v))\right\>}\ff'(\cdot+\al(v))$, and
\item $\<\al'(v),v\>=0$ and $\al'(v)={\|\ff'\|}_\lt^{-2}v'$, if $v\in\mathcal{O}_\ff$.
\end{enumerate}
\end{lemma}
\proof
The proof follows the line of reasoning laid down in Theorem 3.1 in \cite{ribeiro} and Lemma 3.8 in \cite{liutom}.
\fim
Let $\psi$ be as in  Theorem \ref{instability-cond-theo}. Define another vector field $B_\psi$ by
\[
B_\psi(u)=\tau_{\al(u)}\psi'-\frac{\left\<u,\tau_{\al(u)}\psi'\right\>}{\left\<u,\tau_{\al(u)}\ff''\right\>}\tau_{\al(u)}\ff'',
\]
for $u\in\Omega_{\ff,\eps}$. Geometrically, $B_\psi$  can be interpreted as the derivative of the orthogonal component of $\tau_{\al(\cdot)}\psi$ with regard to $\tau_{\al(\cdot)}\ff'$.

\begin{lemma}\label{vector-field}
Let $\psi$ be as in  Theorem \ref{instability-cond-theo}. Then the map $B_\psi:\Omega_{\ff,\eps_0}\to\x$ is $C^1$ with bounded derivative. Moreover,
\begin{enumerate}[(i)]
\item $B_\psi$ commutes with translations,
\item $\<B_\psi(u),u\>=0$, if $u\in \Omega_{\ff,\eps_0}$,
\item $B_\psi(\ff)=\psi'$, if $\<\ff,\psi'\>=0$.
\end{enumerate}
\end{lemma}
\proof
The proof follows the same lines from the proof of Lemma 3.5 in \cite{angulo-1}, Lemma 3.3 in \cite{angulo-2} or Lemma 4.7 in \cite{levaliu}.
\fim

\noindent\textbf{Proof of Theorem \ref{instability-cond-theo}.}\q
First we claim that there exist $\eps_3>0$ and $\si_3>0$ such that for each $u_0\in\Omega_{\ff,\eps_3}$,
\begin{equation}\label{taylor-exp-1}
S(\ff)\leq S(u_0)+\p(u_0)s,
\end{equation}
for some $s\in(-\si_3,\si_3)$, where $\p(u)=\<S'(u),B_\psi(u)\>$.

We consider $u_0\in \in\Omega_{\ff,\eps_0}$, where $\eps_0$ is given in Lemma \ref{implicit}, the initial value problem
\begin{equation}\label{initial-value}
\begin{array}{ll}
\frac{\dd }{\dd s}u(s)=B_\psi(u(s))\\
u(0)=u_0.
\end{array}
\end{equation}
By Lemma \ref{vector-field}, we have that \eqref{initial-value} admits for each $u_0\in\Omega_{\ff,\eps_0}$ a unique maximal solution $u\in C^2((-\si,\si);\Omega_{\ff,\eps_0})$, where $\si\in(0,+\infty]$. Moreover for each $\eps_1<\eps$, there exists $\si_1>0$ such that $\si(u_0)\geq\si_1$, for all $u_0\in \Omega_{\ff,\eps_1}$. Hence we can define for fixed $\eps_1$, $\si_1$, the following dynamical system
\[
\begin{array}{cc}
\uu:(-\si_1,\si)\times \Omega_{\ff,\eps_1}\longrightarrow\Omega_{\ff,\eps_0}\\
\q\q(s,u_0)\mapsto\uu(s)u_0,
\end{array}\]
where $s\to\uu(s)u_0$ is the maximal solution of \eqref{initial-value} with initial data $u_0$. It is also clear from Lemma \ref{vector-field} that $\uu$ is a $C^1-$function, also we have that for each $u_0\in \Omega_{\ff,\eps_1}$, the function $s\to\uu(s)u_0$ is $C^2$ for each $s\in(-\si_1,\si_1)$, and the flow $s\to\uu(s)u_0$ commutes with translations. One can also observe from the relation
\[
\uu(t)\ff=\ff+\int_0^t\tau_{\al(\uu(s)\ff)}\psi'\dd s-\int_0^t\rho(s)\tau_{\al(\uu(s)\ff)}\ff''\dd s
\]
that $\uu(s)\ff\in X_r$, $r>3/2$, for all $s\in(-\si_1,\si_1)$, where
\[
\rho(s)=\frac{\left\<\uu(s)\ff,\tau_{\al(\uu(t)\ff)}\psi'\right\>}{\left\<\uu(t)\ff,\tau_{\al(\uu(t)\ff)}\ff''\right\>}.
\]
Now we get from Taylor's theorem that there is $\varrho\in(0,1)$ such that
\[
S(\uu(s)u_0)=S(u_0)+\p(u_0)s+\frac{1}{2}R(\uu(\varrho s)u_0)s^2,
\]
where $R(u)=\<S''(u)B_\psi,B_\psi(u)\>+\<S''(u),B_\psi'(u)(B_\psi(u))\>$. Since $R$ and $\p$ are continuous, $S'(\ff)=0$ and $R(\ff)<0$, then there exists $\eps_2\in(0,\eps_1]$ and $\si_2\in(0,\si_1]$ such that \eqref{taylor-exp-1} holds for $u_0\in B(\ff,\eps_2)$ and $s\in (-\si_2,\si_2)$. On the other hand, it is straightforward to verify that
\[
\left.P(\uu(s)u_0)\right|_{(u_0,s)=(\ff,0)}=0\q\mbox{and}\q
\frac{\dd}{\dd s}\left.P(\uu(s)u_0)\right|_{(u_0,s)=(\ff,0)}=\<P'(\ff),\psi'\>,
\]
where $P$ is defined in Theorem \ref{ground}. We show that $\<P'(\ff),\psi'\>\neq0$. Otherwise, $\psi'$ would be tangent to $\n$ at $\ff$, is defined in Theorem \ref{ground}. Hence, $\<S''(\ff)\psi',\psi'\>\geq0$, since $\ff$ minimizes $S$ on $\n$ by Theorem \ref{ground}. But this contradicts \eqref{instability-cond}. Therefore, by the implicit function theorem, there exist $\eps_3\in(0,\eps_2)$
and $\si_3\in(0,\si_2)$ such that for all $u_0\in B{\ff,\eps_3}$, there exists a unique $s=s(u_0)\in(-\si_3,\si_3)$ such that $P(\uu(s)u_0)=0$. Then applying \eqref{taylor-exp-1} to $(u_0,s(u_0))\in B{\ff,\eps_3}\times(-\si_3,\si_3)$ and using the fact  $\ff$ minimizes $S$ on $\n$, we have that for $u_0\in B{\ff,\eps_3}$ there exists $s\in(-\si_3,\si_3)$ such that $S(\ff)\leq S(\uu(s)u_0)\leq S(u_0)+\p(u_0)s$. This inequality can be extended to $\Omega_{\ff,\eps_3}$ from from the gauge
invariance.

Since $\uu(s)u_0$ commutes with $\tau_r$, it follows by replacing $u_0$ with $\uu(s)u_0$ in \eqref{taylor-exp-1} and then $\delta=-s$ that
\begin{equation}\label{taylor-exp-2}
S(\ff)\leq S(\uu(\delta)\ff)-\p(\uu(\delta)\ff)\delta,
\end{equation}
for all $\delta\in(-\si_3,\si_3)$. Moreover, using Taylor's theorem again and the fact  $\p(\ff)=0$, it follows that
the map $\delta\mapsto S(\uu(\delta)\ff)$ has a strict local maximum at $\delta=0$. Hence, we obtain
\begin{equation}\label{taylor-exp-3}
S(\uu(\delta)\ff)<S(\ff),\q\q\delta\neq0,\;\delta\in(-\si_4,\si_4),
\end{equation}
where $\si_4\in(0,\si_3]$. Thus it follows from \eqref{taylor-exp-2} that
\begin{equation}\label{taylor-exp-4}
\p(\uu(\delta)\ff)<0,\q\q\delta\in(0,\si_4).
\end{equation}
Let $\delta_j\in(0, \si_4)$ such that $\delta_j\to0$ as $j\to\infty$. Consider the sequences of initial data $u_{0,j}=\uu(\delta_j)\ff$. It is clear to see that $u_{0,j}\in X_s$, $s > 3/2$ for all positive integers $j$ and $u_{0,j}\to\ff$ in $\x$ as $j\to\infty$.

Now we need only verify that the solution $u_j(t)=\uu(t)u_{0,j}$ of \eqref{rgbo} with $u_j(0) = u_{0,j}$ escapes from $\Omega_{\ff,\eps_3}$, for all positive integers $j$ in finite time.  Define
\[
T_j=\sup\{t'>0;\,u_j(t)\in \Omega_{\ff,\eps_3},\;\forall t\in(0,t')\}
\]
and
\[
\ddd=\{u\in\Omega_{\ff,\eps_3};\,S(u)<S(\ff),\;\p(u)<0\}.
\]
Hence it follows from \eqref{taylor-exp-1} that for all $j\in\N$ and $t\in (0,T_j)$, there exists $s=s_j(t)\in(-\si_3,\si_3)$ satisfying $S(\ff)\leq S(u_{0,j})+\p(u_j(t))s$. By \eqref{taylor-exp-3} and \eqref{taylor-exp-4}, $u_{0,j}\in\ddd$; and therefore $u_j(t)\in\ddd$ for all $t\in[0,T_j]$. Indeed, if $\p(u_j(t_0))>0$ for some $t_0\in[0, T_j]$, then the continuity of $\p$ implies that there exists some $t_1\in[0, T_j]$ satisfying $\p(u_j(t_1))=0$, and consequently  $S(\ff)\leq S(u_{0,j})$, which contradicts $u_{0,j}\in\ddd$. Hence,  $\ddd$ is bounded away from zero and
\begin{equation}\label{bound}
-\p(u_j)\geq\frac{S(\ff)-S(u_{0,j})}{\si_3}-\eta_j>0,\q\forall t\in [0,T_j].
\end{equation}
Now suppose that for some $j$, $T_j=+\infty$. Then we define a Liapunov function
\[
A(t)=\int_\rr\psi(x+\al(u_j))u_j(x,t)\dd x,\q t\in[0,T_j].
\]
Then by the Cauchy-Schwarz inequality,
\[
|A(t)|\leq\|\psi\|_\lt\|u_j(t)\|_\lt=\|\psi\|_\lt\|u_{0,j}\|_\lt<\infty,\q\q t\in[0,T_j].
\]
On the other hand, since $\frac{\dd u_j}{\dd t}=-\partial_x E'(u_j)$, then we have
\[
\begin{split}
\frac{\dd A}{\dd t}&=\left\<\al'(u_j(t)),\frac{\dd u_j}{\dd t}\right\>
\<\tau_{\al(u_j(t))}\psi',u_j(t)\>
+\left\<\tau_{\al(u_j(t))},\frac{\dd u_j}{\dd t}\right\>\\
&=\left\<\<\tau_{\al(u_j(t))}\psi',u_j(t)\>\partial_x\al'(u_j(t))+\tau_{\al(u_j(t))}\psi',E'(u_j(t))\right\>\\
&=\left\<B_\psi(u_j(t)),S'(u_j(t))\right\>+c\<B_\psi(u_j(t)),u_j(t)\>
=\p(u_j(t)),
\end{split}
\]
for $t\in[0,T_j]$. Therefore it is deduced from \eqref{bound} that
\[
-\frac{\dd A}{\dd t}\geq\eta_j>0,\q\q\forall t\in[0,T_j].
\]
This contradicts the boundedness of $A(t)$. Consequently $T_j<+\infty$ for all $j$, which means that $u_j$ eventually leaves $\Omega_{\ff,\eps_3}$. This completes the proof.
\fim

\begin{theorem}\label{instability-theo-2} Fix $\beta>0$, $\gamma>0$ and assume there exists a $C^2$ map $c\mapsto\varphi_c\in\g$ for $c<c_*$. If $d''(c)<0$, then $\mathcal{O}_{\ff_c}$ is $\x$-unstable.
\end{theorem}

\proof It suffices to show that there exists a function $\psi$ that satisfies the conditions of Theorem \ref{instability-cond-theo}. Define
\[
\psi(x)=\int_{-\infty}^x\ff_c(y)-\frac{2d'(c)}{d''(c)}\frac{d}{dc}\ff_c(y)\dd y.
\]
Then since $\ff_c\in\x$ and $\frac{d}{dc}\ff_c\in\x$ it follows that $\psi'\in\x$, and thus $\psi\in L^2$. Since $w=\frac{d}{dc}\ff_c$ satisfies the linear equation
    \[
    \beta\h(w_x)-cw-\gamma\partial_x^{-1}w-f'(\ff)w=\ff
    \]
it follows as in the proof of Theorem \ref{regularity} that $w\in H^\infty$ and $\partial_x^{-1}w\in H^\infty$. Hence $\psi'\in X_s$ and $\psi''\in\x$. Now since $d'(c)=-\frac12\<\ff_c,\ff_c\>$, we have
    \[
    \<\psi',\ff_c\>=\<\ff_c,\ff_c\>-\frac{2d'(c)}{d''(c)}\frac12\frac{d}{dc}\<\ff_c,\ff_c\>
    =-2d'(c)+\frac{2d'(c)}{d''(c)}d''(c)=0.
    \]
Next we compute
\[
  \<S''(\ff)\psi',\psi'\>=\<S''(\ff_c)\ff_c,\ff_c\>-\frac{4d'(c)}{d''(c)}\left\<S''(\ff_c)\ff_c,\frac{d}{dc}\ff_c(y)\right\>
  +\frac{4d'(c)^2}{d''(c)^2}\left\<S''(\ff_c)\frac{d}{dc}\ff_c(y),\frac{d}{dc}\ff_c(y)\right\>
\]
For any $\ff\in\g$ we have $S''(\ff)\ff=(p-1)f(\ff)$, so it follows that
\begin{equation}
  \<S''(\ff)\ff,\ff\>=(1-p)K(\ff).
\end{equation}
Since $d(c)=\frac{p-1}{2(p+1)}K(\ff_c)$ we have
    \[
    \left\<S''(\ff)\ff_c,\frac{d}{dc}\ff_c(y)\right\>=\left\<(p-1)f(\ff_c),\frac{d}{dc}\ff_c(y)\right\>
    =-\frac{p-1}{p+1}\frac{d}{dc}K(\ff_c)=-2d'(c).
    \]
Finally, since $S''(\ff_c)\frac{d}{dc}\ff_c(y)=\ff_c$ we have
    \[
    \left\<S''(\ff_c)\frac{d}{dc}\ff_c(y),\frac{d}{dc}\ff_c(y)\right\>
    =\left\<\ff_c,\frac{d}{dc}\ff_c(y)\right\>=-d''(c).
    \]
Altogether this implies
    \[
    \<S''(\ff)\psi',\psi'\>=(1-p)K(\ff_c)+\frac{4(d'(c))^2}{d''(c)}.
    \]
Since $p>1$ and $d''(c)<0$, both terms on the right hand side are negative. Thus $\psi'$ satisfies all of the conditions of Theorem \ref{instability-cond-theo}.
\fim

\section{Applications of the Stability and Instability Theorems}

In this section we apply the stability and instability conditions in Theorems \ref{stability-theo}, \ref{instability-cond-theo} and \ref{instability-theo-2} to determine conditions on $p$, $\beta$, $c$ and $\gamma$ that imply stability or instability. We first apply Theorem \ref{instability-cond-theo} with
$\psi'=\ff+2x\ff'$.

\begin{lemma}\label{quantity-value}
Let $c<c_\ast$ and $\ff\in\g$. Define
\[
\psi(x)=\int_{-\infty}^x\ff(y)+2y\ff'(y)\dd y.
\]
Then $\psi$ satisfies the assumptions of  Theorem \ref{instability-cond-theo} and
\[
\left\<S''(\ff)\psi',\psi'\right\>=\frac{(p-1)(3-p)}{p+1}K(\ff)+12\ga\int_\rr(\dn\ff)^2\dd x.
\]
\end{lemma}
\proof
The first part of the lemma is clear by using the fact $\<\ff,\psi'\>=0$ and Theorems \ref{regularity} and \ref{decay}.
Now we  estimate the quantity $\left\<S''(\ff)\psi',\psi'\right\>$. First by \eqref{rgbo-1}, we note that $S''=\beta\h\partial_x-\gamma\partial_x^{-2}-c-f'(\ff)$. Next, using \eqref{E:functional_identities}, we see that
\begin{equation}\label{quant-1}
\<S''(\ff),\ff\>=(1-p)K(\ff).
\end{equation}
Next using again \eqref{E:functional_identities} an the facts $F'=f$ and $pf(\ff)=f'(\ff)\ff$, it yields that
\begin{equation}\label{quant-2}
\<S''(\ff),x\ff'\>=\int_\rr(p-1)x\ff'f(\ff)\dd x=\frac{p-1}{1+p}K(\ff).
\end{equation}
Finally we show that $\<S''(x\ff'),x\ff'\>=3\ga\int_\rr(\dn\ff)^2\dd x$.

First we observe from \eqref{rgbo-1} and \eqref{E:functional_identities} that
\[\begin{split}
S''(x\ff')&=\beta\h(x\ff')_x-\ga\partial_x^{-2}(x\ff')-cx\ff'-x\ff'f'(\ff)\\
&=
\beta\h\ff'+x\left(\beta\h\ff'-c\ff-\ga\partial_x^{-2}\ff+f(\ff)\right)_x+2\ga\partial_x^{-2}\ff\\
&=\beta\h\ff'+2\ga\partial_x^{-2}\ff=3\ga\partial_x^{-2}\ff+c\ff-f(\ff);
\end{split}\]
and by using \eqref{E:functional_identities} again, we obtain
\begin{equation}\begin{split}\label{quant-3}
\<S''(x\ff'),x\ff'\>&=\int_\rr\left(3\ga\partial_x^{-2}\ff+c\ff-f(\ff)\right)x\ff'\dd x\\
&=\frac{1}{2}\int_\rr\left(9(\dn\ff)^2-c\ff^2\right)\dd x-\frac{1}{p+1}K(\ff)
=3\ga\int_\rr(\dn\ff)^2\dd x.
\end{split}\end{equation}
Therefore we deduce from \eqref{quant-1},\eqref{quant-2} and \eqref{quant-3} that
\[\begin{split}
\<S''(\psi'),\psi'\>&=
\<S''(\ff),\ff\>+4\<S''(\ff),x\ff'\>+\<S''(x\ff'),x\ff'\>\\
&=(1-p)K(\ff)+\frac{4(p-1)}{p+1}K(\ff)+12\ga\int_\rr(\dn\ff)^2\dd x\\
&=\frac{(p-1)(3-p)}{p+1}K(\ff)+12\ga\int_\rr(\dn\ff)^2\dd x.
\end{split}\]
\fim

\begin{theorem}\label{main-instability}
Let $\beta>0$, $\gamma>0$, $c < c_\ast=3(\beta^2\ga/4)^{1/3}$ and $\ff\in\g$. Then the orbit $\mathcal{O}_\ff$ is $\x$-unstable if one the following cases occurs:
 \begin{enumerate}[(i)]
\item $c<0$, $p>3$ and $\ga$ is sufficiently small,
\item $p>5$ and $c<\left(\frac{p-5}{p-1}\right)c_\ast$
\end{enumerate}
\end{theorem}
\proof
By Theorem \ref{instability-cond-theo} and Lemma \ref{quantity-value}, we only need to check condition \eqref{instability-cond} for $\psi$ defined in Lemma \ref{quantity-value}.

First we note that $\lim_{\ga\to0}\ga\int_\rr(\dn\ff)^2\dd x=0$. Indeed, we already know from \eqref{E:functional_identities} that
\[
\ga\int_\rr(\dn\ff)^2\dd x=\int_\rr c\ff^2-\beta(\dx\ff)^2\dd x+(m(\beta,c,\ga))^{\frac{p+1}{p-1}}.
\]
Applying Theorem \ref{weak}, it transpires that
\[
\lim_{\ga\to0}\ga\int_\rr(\dn\ff)^2\dd x=\int_\rr c\phi^2-\beta(\dx\phi)^2\dd x+
(m(\beta,c,0))^{\frac{p+1}{p-1}}
=
-I(\phi;\beta,c,0)+(m(\beta,c,\ga))^{\frac{p+1}{p-1}}=0,
\]
where $\phi$ is a ground state of \eqref{gbo} with $c<0$. Applying Theorem \ref{weak} once more we see that
    \[
    \lim_{\gamma\to0^+}\frac{(p-1)(3-p)}{p+1}K(\ff)+12\ga\int_\rr(\dn\ff)^2\dd x
    =\frac{(p-1)(3-p)}{p+1}m(\beta,c,0)^{\frac{p+1}{p-1}}<0
    \]
since $p>3$. Therefore by Lemma \ref{quantity-value} one has $\<S''(\ff)\psi',\psi'\><0$ for $\gamma>0$ sufficiently small. This which proves (i).

Attention is now given to the proof of (ii). Suppose $p>5$. By Lemma \ref{quantity-value} and equation \eqref{E:functional_identities} we have
    \[
    \<S''(\psi'),\psi'\>=(5-p)K(\varphi)+4c\int_\rr\ff^2\dd x.
    \]
This is clearly negative when $c\leq 0$. Now for $c>0$, a straightforward calculation reveals that
    \[
    \int_\rr\ff^2\dd x\leq\frac1{c_\ast-c}I(\varphi)
    \]
and thus
    \[
    \<S''(\psi'),\psi'\>\leq \left(5-p+\frac{4c}{c_\ast-c}\right)K(\varphi).
    \]
The term on the right hand side is negative when $c<\left(\frac{p-5}{p-1}\right)c_\ast$.
This completes the proof.
\fim

\begin{remark} Notice that as $p\to\infty$, $\left(\frac{p-5}{p-1}\right)c_\ast\to c_\ast$, so the region of instability approaches the entire domain of existence.
\end{remark}

We now investigate what conclusions may be drawn from Theorems \ref{stability-theo} and \ref{instability-theo-2}, which state that stability is determined by the sign of $d''(c)$. Although no explicit formula for $d$ is available, it is possible to determine the behavior of $d''(c)$ for small $\gamma>0$. The following scaling property is the main ingredient in this analysis.

\begin{lemma}\label{L:d_scaling}
Let $\beta>0$, $\gamma>0$ and
$c<c_\ast$. For any $r>0$ and $s>0$ we have
    \[
    d(r\beta,rcs^{-1},rs^{-3}\gamma)=r^{\frac{p+1}{p-1}}s^{\frac{-2}{p-1}}d(\beta,c,\gamma).
    \]
\end{lemma}

\proof The lemma follows from \eqref{E:dIKM_relation} once we show that
    \[
    m(r\beta,rcs^{-1},rs^{-2}\gamma)=rs^{\frac{-2}{p+1}}m(\beta,c,\gamma).
    \]
Let $u\in\x$ with $K(u)>0$. For any $r>0$ we have
    \[
    I(u;r\beta,rc,r\gamma)=rI(u;\beta,c,\gamma),
    \]
so $m(r\beta,rc,r\gamma)=rm(\beta,c,\gamma)$. Next let
$v(x)=u(sx)$ for $s>0$. Then
    \[
    I(v;\beta,c,\gamma)=I(u;\beta,cs^{-1},s^{-3}\gamma)
    \qquad K(v)=\frac1s K(u)
    \]
so
    \[
    \frac{I(v;\beta,c,\gamma)}{K(v)^{\frac{2}{p+1}}}
    =s^{\frac{2}{p+1}}\frac{I(u;\beta,cs^{-1},s^{-3}\gamma)}{K(u)^{\frac{2}{p+1}}}
    \]
and consequently
    \[
    m(\beta,cs^{-1},s^{-3}\gamma)=s^{\frac{-2}{p+1}}m(\beta,c,\gamma).
    \]
\fim
Setting $r=2/\beta$ and $s^3=2\gamma/\beta$ gives
    \begin{equation}\label{E:d_scaling1}
    d\left(2,c\left(\frac4{\gamma\beta^2}\right)^{1/3},1\right)=\left(\frac{2\gamma}{\beta}\right)^{\frac{-2}{3(p-1)}}
    \left(\frac{2}{\beta}\right)^{\frac{p+1}{p-1}}d\left(\beta,c,\gamma\right)
    \end{equation}
Hence for any constant $k$, the values of $d$ along the surface $c^3=k\gamma\beta^2$ are determined by the value of $d$ at any single point on that surface. Next, setting $r=1$ and $s=-c/3$ gives
    \begin{equation}\label{E:d_scaling2}
    d(\beta,-3,\gamma(-3/c)^3)=\left(-c/3\right)^{\frac{-2}{p-1}}d(\beta,c,\gamma).
    \end{equation}
or equivalently
    \[
    d(\beta,c,\gamma)=\left(-c/3\right)^{\frac{2}{p-1}}d(\beta,-3,\gamma(-3/c)^3).
    \]
Next we set $q=\frac{2}{p-1}$ and assume that $d$ is twice differentiable. Then differentiating with respect to $c$ gives
\[
d_c(\beta,c,\gamma)=\left(-\frac13q\left(-c/3\right)^{q-1}d+\gamma\left(-c/3\right)^{q-4}d_\gamma
\right)\bigg|_{(\beta,-3,-27\gamma/c^3)}
\]
and
    \begin{equation}\label{d_derivatives}
    d_{cc}(\beta,c,\gamma)=\left(\frac19q(q-1)\left(-c/3\right)^{q-2}d
    -\frac13\gamma(2q-4)\left(-c/3\right)^{q-5}d_\gamma
    +\gamma^2\left(-c/3\right)^{q-8}d_{\gamma\gamma}\right)\bigg|_{(\beta,-3,-27\gamma/c^3)}.
    \end{equation}

\begin{theorem}\label{stability-theo-2} Assume $d$ is twice differentiable on the domain $c<c_\ast$.
\begin{enumerate}[(i)]
\item
Fix $1<p<3$, $\beta>0$ and $c<0$. Then there exist $\gamma_k\to0^+$ such that $d_{cc}(\beta,c,\gamma_k)>0$.
\item Fix $p>3$, $\beta>0$ and $c<0$. Then there exist $\gamma_k\to0^+$ such that $d_{cc}(\beta,c,\gamma_k)<0$.
\end{enumerate}
\end{theorem}

\proof First observe that
    \[
    \lim_{\gamma\to0^+}\frac19q(q-1)\left(-c/3\right)^{q-2}d(\beta,-3,-27\gamma/c^3)
    =\frac19\frac{2(3-p)}{(p-1)^2}\left(-c/3\right)^{(4-2p)/(p-1)}\frac{p-1}{2(p+1)}m(\beta,-3,0)^{\frac{p+1}{p-1}}.
    \]
This is positive when $1<p<3$ and negative when $p>3$. As shown in the proof of Theorem \ref{main-instability}, the term
    \[
    \gamma d_\gamma=\gamma\int_\rr(\dn\ff)^2\dd x
    \]
vanishes as $\gamma$ approaches zero. It therefore remains to show that the term $\gamma^2 d_{\gamma\gamma}$ vanishes as well. To do so, define
    \[
    g(\gamma)=
    \left\{\begin{array}{cc}
    \gamma^2d_\gamma&\gamma>0\\
    0&\gamma=0
    \end{array}
    \right.
    \]
Then since $\gamma d_\gamma\to0$ as $\gamma\to0^+$, $g$ defines a continuous function for $\gamma\geq0$.
Furthermore by the assumption that $d$ is differentiable, it follows that $g$ is differentiable for $\gamma>0$. By the Mean Value Theorem, for each integer $k$ there exists $\gamma_k\in(0,1/k)$ such that $g(1/k)-g(0)=\frac1kg'(\gamma_k)$, and thus
    \[
    g'(\gamma_k)=k g\left(\frac1k\right)=\frac1kd_\gamma\left(\frac1k\right)\to0
    \]
as $k\to\infty$. Now
    \[
    g'(\gamma_k)=2\gamma_k d_\gamma(\gamma_k)+\gamma_k^2d_{\gamma\gamma}(\gamma_k),
    \]
so we have
    \[
    \lim_{k\to\infty}\gamma_k^2d_{\gamma\gamma}(\gamma_k)=\lim_{k\to\infty} g'(\gamma_k)-2\gamma_k d_\gamma(\gamma_k)=0.
    \]
\fim

We next consider the behavior of $d$ for $c$ near $c_\ast=3(\beta^2\ga/4)^{1/3}$. Using appropriately chosen trial functions, we obtain upper bounds on $d$ as $c$ approaches $c_\ast$ for the nonlinearities $f(u)=|u|^p$ and $f(u)=-|u|^{p-1}u$. In both cases, for any $p\geq2$, these bounds imply that $d(c)\to0$ as $c\to c_\ast$. In the case of the odd nonlinearity $f(u)=-|u|^{p-1}u$, the bound implies that $d$ is convex (and hence $\g$ is stable) for $c$ near $c_\ast$.

Our choice of trial function is $u=w_x$, where $w(x)=\ee^{-a|x|}\sin(bx)$ for appropriately chosen $a>0$ and $b\neq0$. It is clear that $u\in\x$, and we have
    \[
    I(u)=\int_\rr (\beta|\xi|^3-c|\xi|^2+\gamma)|\hat{w}|^2\,\dd\xi.
    \]
Since
    \[
    \hat{w}(\xi)=\frac{-4iab\xi}{(\xi^2-(a^2+b^2))^2+4a^2\xi^2},
    \]
we have
    \[
    I(u)=32a^2b^2\int_0^\infty \frac{\beta\xi^5-c\xi^4+\gamma\xi^2}{((\xi^2-(a^2+b^2))^2+4a^2\xi^2)^2}\,\dd\xi.
    \]
This integral may be evaluated explicitly using Maple to obtain
    \begin{equation}\label{E:I_formula}
    I(u)=\frac1{ab(a^2+b^2)}\left(2\beta\arctan\left(\frac{b}{a}\right)\left(a^2+b^2\right)^3
    +\pi(\gamma b^3-cb^5-ca^2b^3)+\beta(2b^5a-2a^5b))\right).
    \end{equation}
For $c<c_\ast$, the cubic $\beta r^3-cr^2+\gamma$ has one real root and two complex roots. Let $b\pm ai$ denote the complex roots. Then by the cubic formula, we have
    \[
    a=\frac{\sqrt{3}}{6\beta}\left(\frac{D}{2}-\frac{2c^2}{D}\right)
    \]
and
    \[
    b=-\frac{D}{12\beta}-\frac{c^2}{3\beta D}+\frac{c}{3\beta},
    \]
where
    \[
    D=(8c^3-108\gamma\beta^2+12\beta\sqrt{3\gamma(27\gamma\beta^2-4c^3)})^{1/3}.
    \]
As $c\to c_\ast=3(\beta^2\gamma/4)^{1/3}$, we have $D\to -2c_\ast$ and thus
\begin{align*}
  \lim_{c\to c_\ast}a&=0\\
  \lim_{c\to c_\ast}b&=\frac{2c_\ast}{3\beta}.
\end{align*}
Moreover,
    \[
    D+2c_\ast=\frac{D^3+8c_\ast^3}{D^2-2c_\ast D+4c_\ast^2}
    =8(c^3-c_\ast^3)+12\beta\sqrt{12\gamma(c_\ast^3-c^3)}=O(\sqrt{c_\ast-c}),
    \]
    \begin{align*}
    a&=\frac{\sqrt{3}(D^2-4c^2)}{12\beta D}=\frac{\sqrt{3}(D-2c)(D^3+8c^3)}{12\beta D(D^2-2cD+4c^2)}\\
    &=\frac{\sqrt{3}(D-2c)}{12\beta D(D^2-2cD+4c^2)}\left(16(c^3-c_\ast^3)+12\beta\sqrt{12\gamma(c_\ast^3-c^3)}\right)\\
    &=O(\sqrt{c_\ast-c})
    \end{align*}
and
    \begin{align*}
    b-\frac{2c_\ast}{3\beta}&=\frac{c-c_\ast}{3\beta}-\frac{D^2+4c^2+4c_\ast D}{12\beta D}\\
    &=\frac{c-c_\ast}{3\beta}-\frac{4(c^2-c_\ast^2)+(D+2c_\ast)^2}{12\beta D}\\
    &=O(c-c_\ast)
    \end{align*}
as $c\to c_\ast$.
\begin{lemma}
Suppose that $d(c)$ is differentiable for $c<c_\ast$. Then for $c<c_\ast$ it holds that
\begin{equation}\label{lowerbound}
d(c)\geq d(0)\left(1-\frac{c}{c_\ast}\right)^{\frac{p+1}{p-1}}.
\end{equation}
\end{lemma}
\proof
By \eqref{best-estimate}, \eqref{E:dIKM_relation} and Lemma \ref{L:d_derivatives}, it follows that
\[
d(c)=\frac{p-1}{2(p+1)}I(\ff)\geq\frac{p-1}{2(p+1)}(c_\ast-c)\int_\rr\ff^2\dd x
=-\frac{p-1}{p+1}(c_\ast-c)d'(c).
\]
Hence, we obtain that
\[
\frac{d'(c)}{d(c)}\geq\frac{p+1}{(p-1)(c-c_\ast)},
\]
and therefore \eqref{lowerbound} follows.
\fim
\begin{lemma} For $u$, $a$ and $b$ as chosen above, we have
    \[
    I(u)=O(\sqrt{c_\ast-c})
    \]
as $c\to c_\ast$.
\end{lemma}

\proof
Since $a=O(\sqrt{c_\ast-c})$ and $b=O(1)$ as $c\to c_\ast$, it suffices to show that the term in parentheses in expression \eqref{E:I_formula} is $O(c_\ast-c)$. Using the expansion
    \[
    \arctan\left(\frac1{x}\right)=\frac{\pi}{2}-x+O(x^2)
    \]
which holds for small $x>0$, we have
    \[
    2\beta\arctan\left(\frac{b}{a}\right)=\beta\pi-2\beta\frac{a}{b}+O(a^2/b^2)
    \]
and thus
    \[
    2\beta\arctan\left(\frac{b}{a}\right)(a^2+b^2)^3=\beta\pi b^6-2\beta ab^5+O(a^2).
    \]
Combining this with the other two terms in equation \eqref{E:I_formula} we are left with
    \[
    \pi b^3(\beta b^3-cb^2+\gamma)+O(a^2)=\pi b^3(\beta b^3-cb^2+\gamma)+O(c_\ast-c).
    \]
Finally, since $b=\frac{2c_\ast}{3\beta}+O(c_\ast-c)$, it follows that
    \begin{align*}
    \beta b^3-cb^2+\gamma&=\beta\left(\frac{2c_\ast}{3\beta}\right)^3
    -c\left(\frac{2c_\ast}{3\beta}\right)^2+\gamma+O(c_\ast-c)\\
    &=(c_\ast-c)\left(\frac{2c_\ast}{3\beta}\right)^2+O(c_\ast-c)\\
    &=O(c_\ast-c).
    \end{align*}
\fim
This bound on $I(u)$, together with a {\it lower} bound on $K(u)$, leads to an upper bound on $m(\beta,c,\gamma)$. The lower bound on $K(u)$ depends on the nonlinear term $f(u)$. For even nonlinearities we have the following bound.

\begin{lemma} Suppose $f(u)=\pm|u|^{p}$. Fix $\beta>0$ and $\gamma>0$. Then
    \[
    d(c)=O(\sqrt{c_\ast-c})
    \]
as $c$ approaches $c_\ast$.
\end{lemma}

\proof It suffices to prove that $K(u)\geq C\sqrt{c_\ast-c}$ for some constant $C$ independent of $c$. For then
    \[
    m(\beta,c,\gamma)\leq\frac{I(u)}{K(u)^{\frac2{p+1}}}\leq \frac{C(c_\ast-c)^{1/2}}{(c_\ast-c)^{\frac1{p+1}}}
    =O(c_\ast-c)^{\frac{p-1}{2(p+1)}}
    \]
and it follows from \eqref{E:dIKM_relation} that
    \[
    d(c)=\frac{p-1}{2(p+1)}m(\beta,c,\gamma)^{\frac{p+1}{p-1}}=O(c_\ast-c)^{\frac12}.
    \]
To obtain the lower bound on $K(u)$, first write
    \[
    K(u)=\int_\rr |u|^pu\,\dd x=2\int_0^\infty \ee^{-a(p+1)x}|b\cos(bx)-a\sin(bx)|^p(b\cos(bx)-a\sin(bx))\,\dd x.
    \]
Rewriting $b\cos(bx)-a\sin(bx)=\sqrt{a^2+b^2}\cos(bx+\phi)$ where $\phi=\arctan(a/b)$ this becomes
    \[
    2(a^2+b^2)^{\frac{p+1}2}\int_0^\infty \ee^{-a(p+1)x}|\cos(bx+\phi)|^p\cos(bx+\phi)\,\dd x,
    \]
and after the change of variable $y=bx+\phi$ this becomes
    \[
    \frac{2\ee^{a(p+1)\phi/b}}{b}(a^2+b^2)^{\frac{p+1}2}\int_\phi^\infty \ee^{-a(p+1)y/b}|\cos(y)|^p\cos(y)\,\dd y.
    \]
As $c$ approaches $c_\ast$ the term outside the integral approaches $2(\gamma/\beta)^{p/4}>0$, so we will henceforth ignore this term. We now break up the integral as
    \[
    \int_\phi^{0}\ee^{-a(p+1)y/b}|\cos(y)|^p\cos(y)\,\dd y+\sum_{k=0}^\infty\int_{k\pi}^{(k+1)\pi}\ee^{-a(p+1)y/b}|\cos(y)|^p\cos(y)\,\dd y.
    \]
The first term is negative, but bounded below by
    \[
    -\phi=-\arctan(a/b)\geq-a/b.
    \]
In each term of the summation we make the change of variable $z=y-k\pi$ to obtain
    \[
    \sum_{k=0}^\infty\int_0^\pi \ee^{-a(p+1)(z+k\pi)/b}|\cos(z)|^p(-1)^k\cos(z)\,\dd z
    \]
which, after summing the geometric series, can be rewritten as
    \[
    \frac1{1+\ee^{-a(p+1)\pi/b}}\int_0^\pi \ee^{-a(p+1)z/b}|\cos(z)|^p\cos(z)\,\dd z.
    \]
The remaining integral we rewrite as
    \[
    \int_0^{\pi/2} \ee^{-a(p+1)z/b}|\cos(z)|^p\cos(z)\,\dd z+\int_{\pi/2}^\pi \ee^{-a(p+1)z/b}|\cos(z)|^p\cos(z)\,\dd z
    \]
and make the change of variable $y=\pi-z$ in the second integral to obtain
    \[
    \int_0^{\pi/2}\ee^{-a(p+1)z/b}|\cos(z)|^p\cos(z)\,\dd z+\int_{\pi/2}^0 \ee^{-a(p+1)(\pi-y)/b}|\cos(y)|^p\cos(y)\,\dd y.
    \]
Combining these, we have
    \[
    \int_0^{\pi/2}(\ee^{-a(p+1)z/b}-\ee^{-a(p+1)(\pi-z)/b})\cos(z)^{p+1}\,\dd z.
    \]
Since
    \[
    \lim_{a\to 0}\frac{\ee^{-a(p+1)z/b}-\ee^{-a(p+1)(\pi-z)/b}}{a}=\frac{p+1}{b}\cdot(\pi-2z)
    \]
uniformly in $x$ in $[0,\pi/2]$ the integral approaches
    \[
    \frac{a}{b}\int_0^{\pi/2}(p+1)(\pi-2z)|\cos(z)|^p\cos(z)\,\dd z
    \]
as $c\to c_\ast$. Since
    \begin{align*}
    \int_0^{\pi/2}(p+1)(\pi-2z)\cos(z)^{p+1}\,\dd z&=2(p+1)\int_0^{\pi/2}x\sin(x)^{p+1}\,\dd x\\
    &\geq 2(p+1)\int_0^{\pi/2}x(2x/\pi)^{p+1}\,\dd x\\
    &=\frac{(p+1)\pi^2}{2(p+3)}\\
    &>\frac{\pi^2}{4}
    \end{align*}
for all $p>1$, it follows that as $c\to c_\ast$ we have
    \[
    \frac1{1+\ee^{-a(p+1)\pi/b}}\int_0^\pi \ee^{-a(p+1)z/b}|\cos(z)|^p\cos(z)\,\dd z\geq \frac12\cdot\frac14\pi^2\cdot\frac{a}{b},
    \]
and therefore
    \[
    \int_\phi^\infty \ee^{-a(p+1)y/b}|\cos(y)|^p\cos(y)\,\dd y\geq \frac12\left(\frac14\pi^2-2\right)\frac{a}{b},
    \]
which implies that
    \[
    K(u)\geq O(a)=O(\sqrt{c_\ast-c})
    \]
as desired.
\fim

While the bound in the previous lemma shows that $d\to0$ as $c\to c_\ast$, unfortunately it does not provide any information about the sign of $d''(c)$. For the odd nonlinearity $f(u)=-|u|^{p-1}u$, however, the integrand of the functional $K$ is nonnegative, and we have the following stronger bound.

\begin{lemma} Suppose $f(u)=-|u|^{p-1}u$. Fix $\beta>0$ and $\gamma>0$. Then
    \[
    d(c)=O(c_\ast-c)^{\frac{p+3}{2(p-1)}}
    \]
as $c$ approaches $c_\ast$.
\end{lemma}

\proof It suffices to prove that $K(u)\geq C(c_\ast-c)^{-1/2}$ for some constant $C$ independent of $c$. For then
    \[
    m(\beta,c,\gamma)\leq\frac{I(u)}{K(u)^{\frac2{p+1}}}\leq \frac{C(c_\ast-c)^{1/2}}{(c_\ast-c)^{-\frac1{p+1}}}
    =O(c_\ast-c)^{\frac{p+3}{2(p+1)}}
    \]
and the lemma follows from \eqref{E:dIKM_relation}. Now, using the calculations from the previous lemma, we have
    \begin{align*}
    K(u)&=\int_\rr |u|^{p+1}\,\dd x=\frac{2\ee^{a(p+1)\phi/b}}{b}(a^2+b^2)^{\frac{p+1}2}\int_\phi^\infty \ee^{-a(p+1)y/b}|\cos(y)|^{p+1}\,\dd y\\
    &\geq \frac{2\ee^{a(p+1)\phi/b}}{b}(a^2+b^2)^{\frac{p+1}2}\int_{\pi/2}^\infty \ee^{-a(p+1)y/b}|\cos(y)|^{p+1}\,\dd y.
    \end{align*}
Writing the integral as
    \[
    \sum_{k=1}^\infty \int_{(k-\frac12)\pi}^{(k+\frac12)\pi} \ee^{-a(p+1)y/b}|\cos(y)|^{p+1}\,\dd y,
    \]
and making the change of variable $z=y-k\pi$, this becomes
    \[
    \sum_{k=1}^\infty \ee^{-a(p+1)\pi k/b}\int_{-\pi/2}^{\pi/2} \ee^{-a(p+1)z/b}\cos(z)^{p+1}\,\dd z
    =\frac{\ee^{a(p+1)\pi/b}}{\ee^{a(p+1)\pi/b}-1}\int_{-\pi/2}^{\pi/2} \ee^{-a(p+1)z/b}\cos(z)^{p+1}\,\dd z.
    \]
For small $a$ this is approximately
    \[
    \frac{b}{a(p+1)\pi}\int_{-\pi/2}^{\pi/2}\cos(z)^{p+1}\,\dd z\geq C'a^{-1}=O(c_\ast-c)^{-1/2}.
    \]
\fim

\begin{theorem}\label{T:d_near_c_star_odd} Suppose $f(u)=-|u|^{p-1}u$ where $1<p<5$. Fix $\beta>0$ and $\gamma>0$. Then there exist $c$ arbitrarily close to $c_\ast$ for which $\g$ is $\x$-stable.
\end{theorem}

\proof For $1<p<5$ the function $(c_\ast-c)^{\frac{p+3}{2(p-1)}}$ is convex and vanishes at $c=c_\ast$. Since $d$ is positive and is bounded above by a multiple of this convex function, its second derivative must be positive at points $c$ arbitrarily close to $c_\ast$.
\fim

\section{Numerical Studies}

In this section we present numerical results which illustrate the behavior of the solitary waves as the parameters $c$ and $\gamma$ are varied, and provide insight into the nature of the function $d(c)$ whose concavity determines the stability of the solitary waves. To obtain the numerical approximations we use a spectral method due to Petviashvili. First observe that the solitary wave equation \eqref{rgbo-1} may be written
    \[
\beta\h\ff_{xxx}-c\ff_{xx}+f(\ff)_{xx}=\ga\ff.
    \]
Writing $\psi_{xx}=\ff$ this becomes
    \[
    -\beta\h\psi_{xxx}+c\psi_{xx}+\gamma\psi=f(\ff)
    \]
so taking the Fourier transform yields
    \[
    (\beta|\xi|^3-c\xi^2+\gamma)\widehat{\psi}=\widehat{f(\ff)}.
    \]
Thus a natural iterative scheme is the following:
\begin{align*}
  \widehat{\psi_{n+1}}&=\frac{\widehat{f(\ff_n)}}{\beta|\xi|^3-c\xi^2+\gamma}\\
  \ff_{n+1}&=(\psi_n)_{xx}.
\end{align*}
Unfortunately, the algorithm has poor convergence properties. However, the algorithm
\begin{align*}
  \widehat{\psi_{n+1}}&=M_n^\alpha\frac{\widehat{f(\ff_n)}}{\beta|\xi|^3-c\xi^2+\gamma}\\
  \ff_{n+1}&=(\psi_n)_{xx},
\end{align*}
with stabilizing factor $M_n$ defined by
    \[
    M_n=\frac{\int(\beta|\xi|^3-c\xi^2+\gamma)|\hat\psi_n|^2\,\dd\xi}{\int\overline{\hat\psi_n}\widehat{f(\ff_n)}\,\dd\xi}
    \]
has much better convergence properties. It was shown in \cite{B:pelinovsky_stepanyants} that this algorithm converges for $1<\alpha<(p+1)/(p-1)$ and the rate of converges is fastest when $\alpha=\alpha^*=p/(p-1)$. This algorithm was implemented in MATLAB using a large spatial domain to compute the solitary waves for a range of parameter values $(\beta,c,\gamma)$. Figures \ref{F:c_limit} and \ref{F:gamma_limit} show several numerically computed solitary waves for the nonlinearity $f(u)=u^2$. Figure \ref{F:c_limit} illustrates the oscillatory tails that develop as $c$ approaches $c_\ast$, while Figure \ref{F:gamma_limit} illustrates the convergence to the exact solitary wave solution of the Benjamin-Ono equation as $\gamma$ approaches zero.

\begin{figure}[ht!]
\begin{center}
  \scalebox{0.5}{\includegraphics{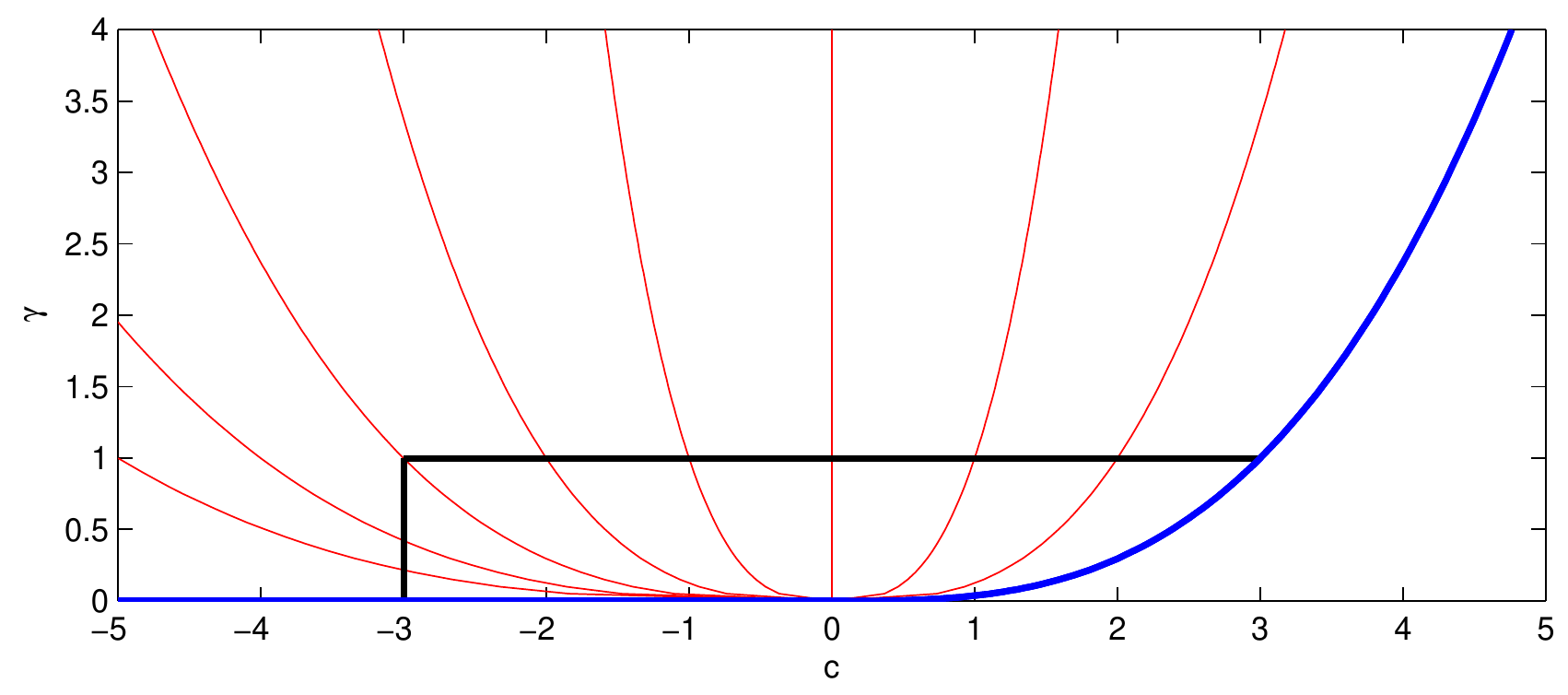}}
\end{center}
\caption{For $\beta=2$, the domain of $d$ is $\{(c,\gamma):\gamma>0, c^3<27\gamma\}$. The numerical computations were performed along the segments $\{-3\leq c<3,\gamma=1\}$ and $\{c=-3,0<\gamma\leq1\}$. Every curve of the form $c^3=k\gamma$ within the domain of $d$ passes through one of these segments.}\label{F:d_domain}
\end{figure}

Once a solitary wave $\varphi\in\g$ is computed, the values of $d(\beta,c,\gamma)$, $d_c(\beta,c,\gamma)$ and $d_\gamma(\beta,c,\gamma)$ are found by using relation \eqref{E:dIKM_relation} and Lemma \ref{L:d_derivatives}. The domain of $d(\beta,c,\gamma)$ is the region $\{(\beta,c,\gamma):\beta>0,\gamma>0,c^3<27\beta^2\gamma/4\}$, shown in Figure \ref{F:d_domain}. By the scaling relation \eqref{E:d_scaling1}, it suffices to compute $d(\beta,c,\gamma)$ at a single point $(\beta,c,\gamma)$ on each surface of the form $c^3=k\gamma\beta^2/4$ for $k<27$. The segments $S_1=\{\beta=2,\gamma=1,-3\leq c<3\}$ and $S_2=\{\beta=2,c=-3,0<\gamma\leq 1\}$ cross all of these surfaces. Along the segment $S_1$, $d_{cc}$ is computed numerically using the computed values of $d_c$, while along $S_2$, relation \eqref{d_derivatives} is used to compute $d_{cc}$ in terms of the numerically values of $d$, $d_\gamma$ and $d_{\gamma\gamma}$.

\begin{figure}
  \begin{center}
    \scalebox{0.75}{\includegraphics{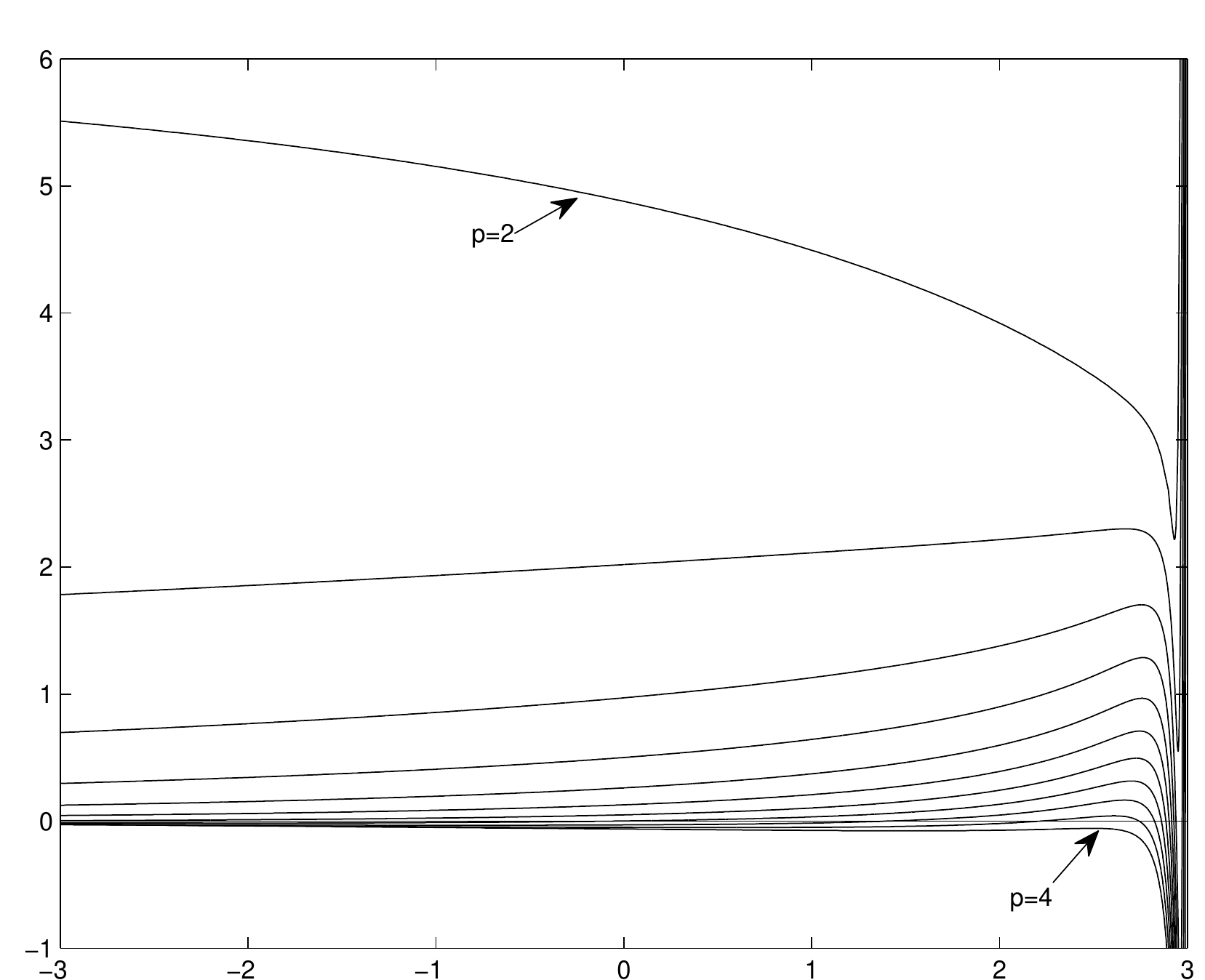}}
    \vskip 5pt
    \scalebox{0.75}{\includegraphics{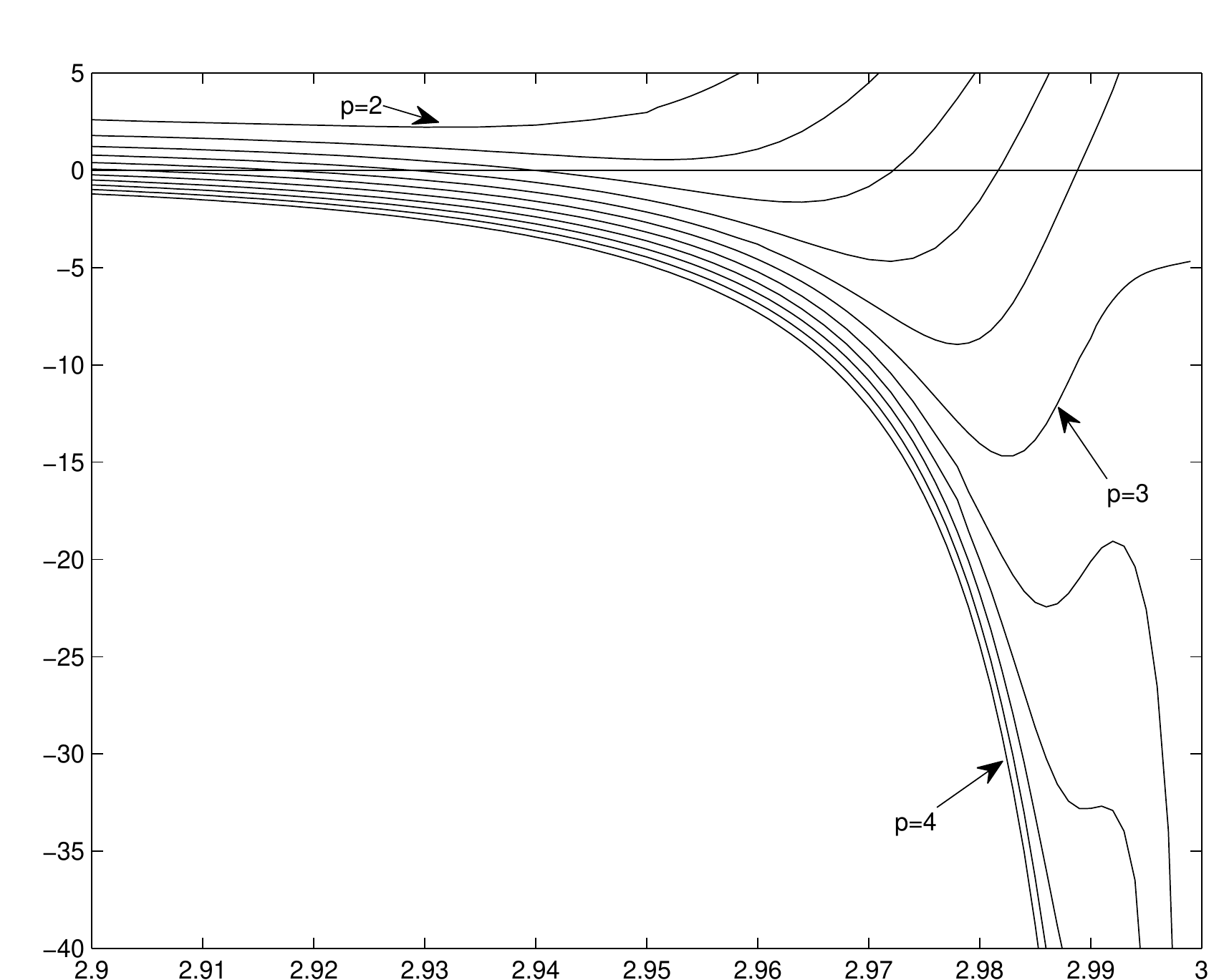}}
  \end{center}
  \caption{Plots of $d_{cc}$ for $f(u)=|u|^p$ with $\beta=2$, $\gamma=1$, $-3\leq c<3$ and $p=2,2.2,2.4,\ldots,4$. The second plot is a blowup of the first, illustrating the behavior for $c$ near $c_\ast=3$}\label{F:d1_odd}
\end{figure}

\begin{figure}
  \begin{center}
    \scalebox{0.45}{\includegraphics{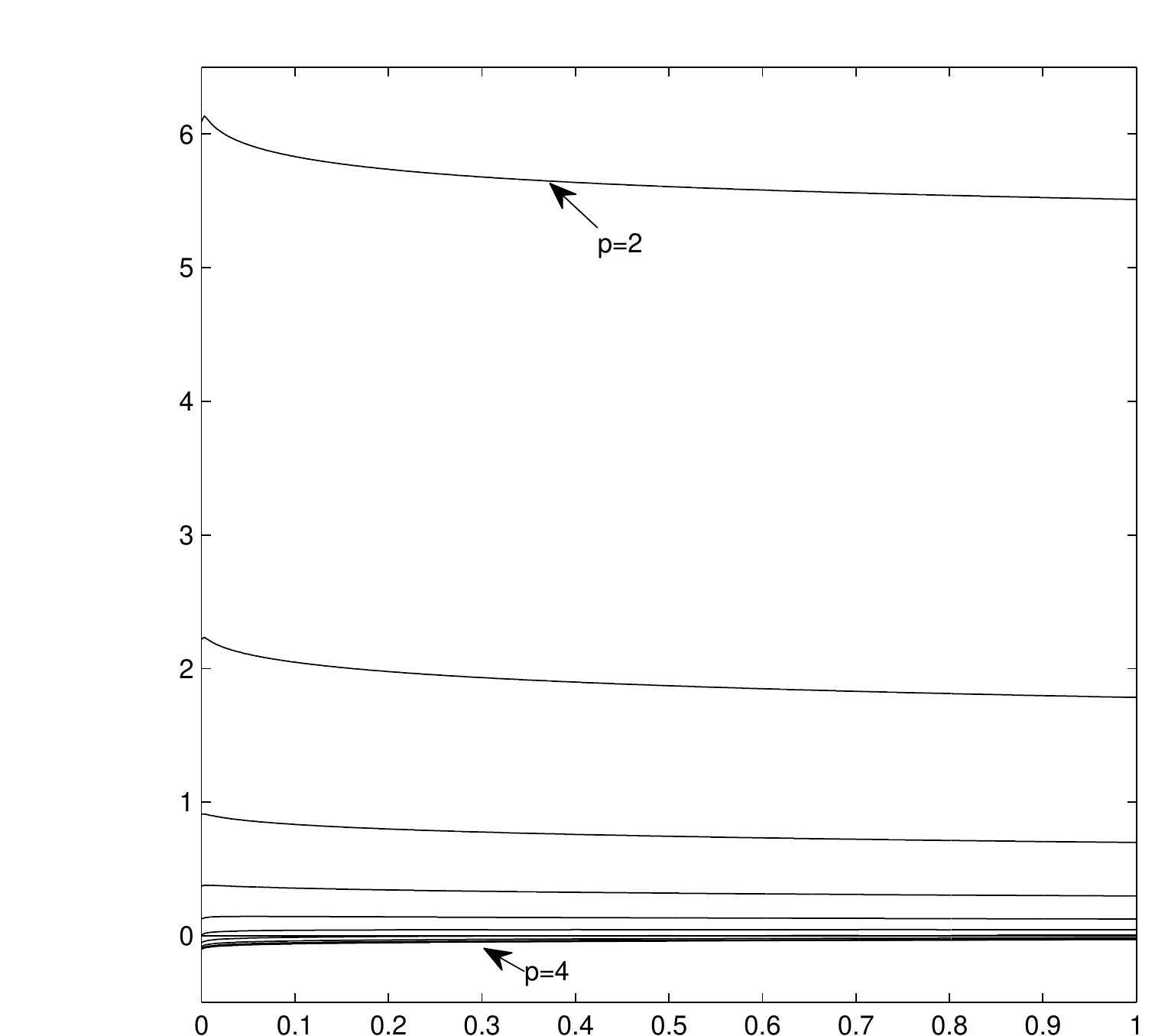}}
    \quad
    \scalebox{0.45}{\includegraphics{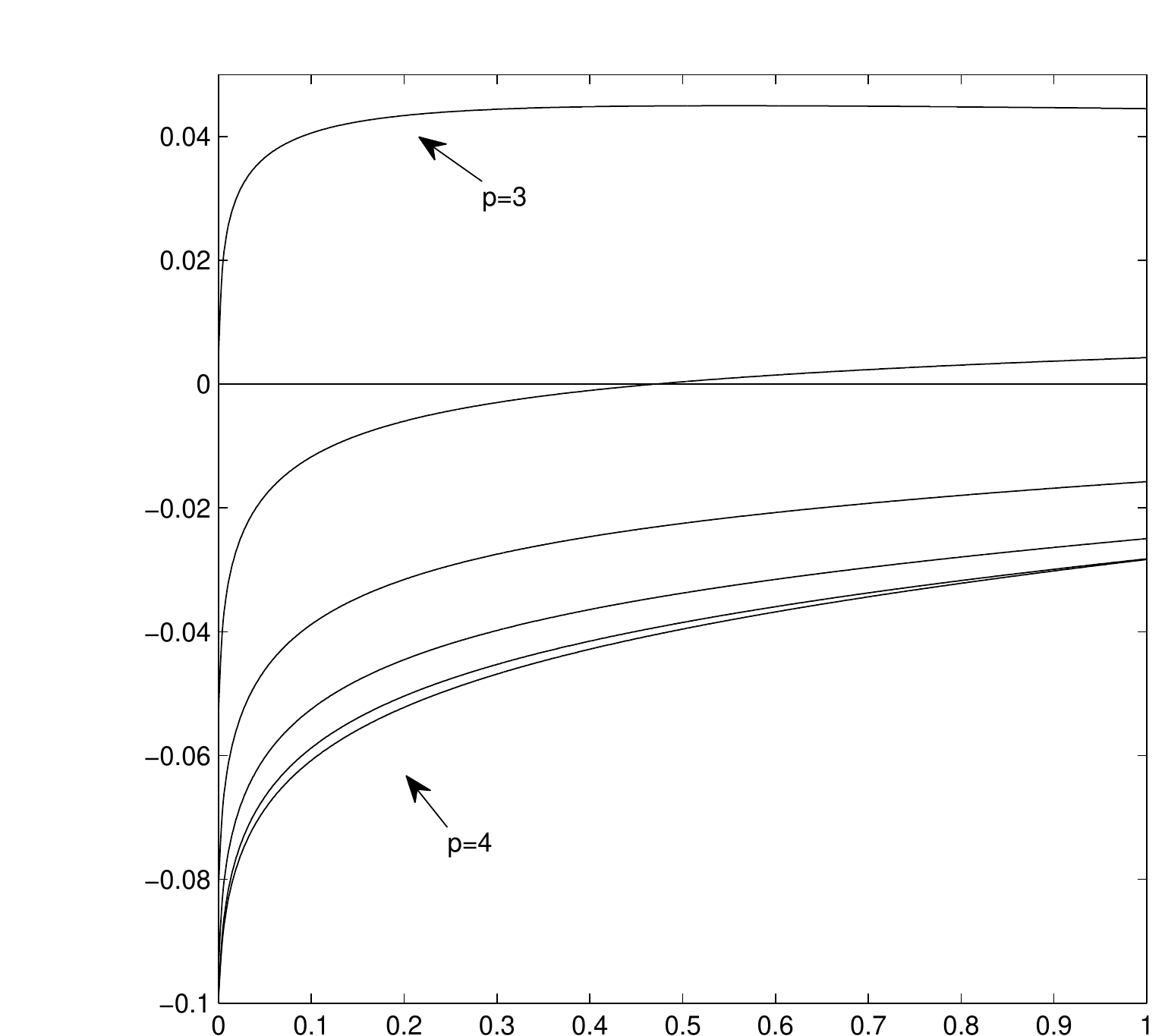}}
 \end{center}
\caption{Plots of $d_{cc}$ for $f(u)=|u|^p$ with $\beta=2$, $c=-3$, $0<\gamma\leq1$ and $p=2,2.2,2.4,\ldots,4$. The second plot is a blowup of the first, and better illustrates the plots for $3\leq p\leq 4$.}\label{F:d2_odd}
\end{figure}

These computations were performed for two families of nonlinearities, even nonlinearities of the form $f(u)=|u|^p$ and odd nonlinearities of the form $f(u)=-|u|^{p-1}u$. The results for the even nonlinearity $f(u)=|u|^p$ are shown in Figures \ref{F:d1_odd} and \ref{F:d2_odd} and summarized in Table \ref{T:d_odd}. For $p=2$ and $p=2.2$ we have $d_{cc}>0$ for all $c<c_\ast$. However, when $p=2.4$ there is a small interval of speeds for which $d_{cc}<0$. As $p$ increases this interval grows, and when $p=4$ we have $d_{cc}<0$ for all $c<c_\ast$. The behavior for small $\gamma>0$ agrees with the results of Theorems \ref{main-instability} and \ref{stability-theo-2} in that when $p<3$ we have $d_{cc}>0$ for small $\gamma$ and when $p>3$ we have $d_{cc}<0$ for small $\gamma$. We note that in the case $p=3$, to which these theorems do not apply, we have $d_{cc}>0$ for small $\gamma$.  The behavior for $c$ near $c_\ast$ is rather interesting. It appears that, for $p<3$, $d_{cc}\to+\infty$ as $c\to c_\ast$, while for $p>3$, $d_{cc}\to-\infty$ as $c\to c_\ast$. When $p=3$, $d_{cc}$ appears to approach some finite negative value.

\begin{table}[h!b!p!]
\caption{Sign of $d_{cc}$ for $f(u)=|u|^p$.}
\begin{center}
\begin{tabular}{|l|c|c|}
\hline
  $p$ & Regions where $d_{cc}>0$.\\
  \hline\hline
  2 & $c<c_\ast$ \\
  \hline
  2.2 & $c<c_\ast$ \\
  \hline
  2.4 & $c<0.980c_\ast$ and $c>0.991c_\ast$\\
  \hline
  2.6 & $c<0.976c_\ast$ and $c>0.994c_\ast$\\
  \hline
  2.8 & $c<0.972c_\ast$ and $c>0.996c_\ast$\\
  \hline
  3 & $c<0.968c_\ast$ \\
  \hline
  3.2 & $-1.287c_\ast<c<0.962c_\ast$ \\
  \hline
  3.4 & $-0.023c_\ast<c<0.954c_\ast$ \\
  \hline
  3.6 & $0.465c_\ast<c<0.942c_\ast$ \\
  \hline
  3.8 & $0.738c_\ast<c<0.915c_\ast$ \\
  \hline
  4 & empty \\
  \hline
  \end{tabular}\label{T:d_odd}
  \end{center}
\end{table}

\begin{figure}
 \begin{center}
   \scalebox{0.75}{\includegraphics{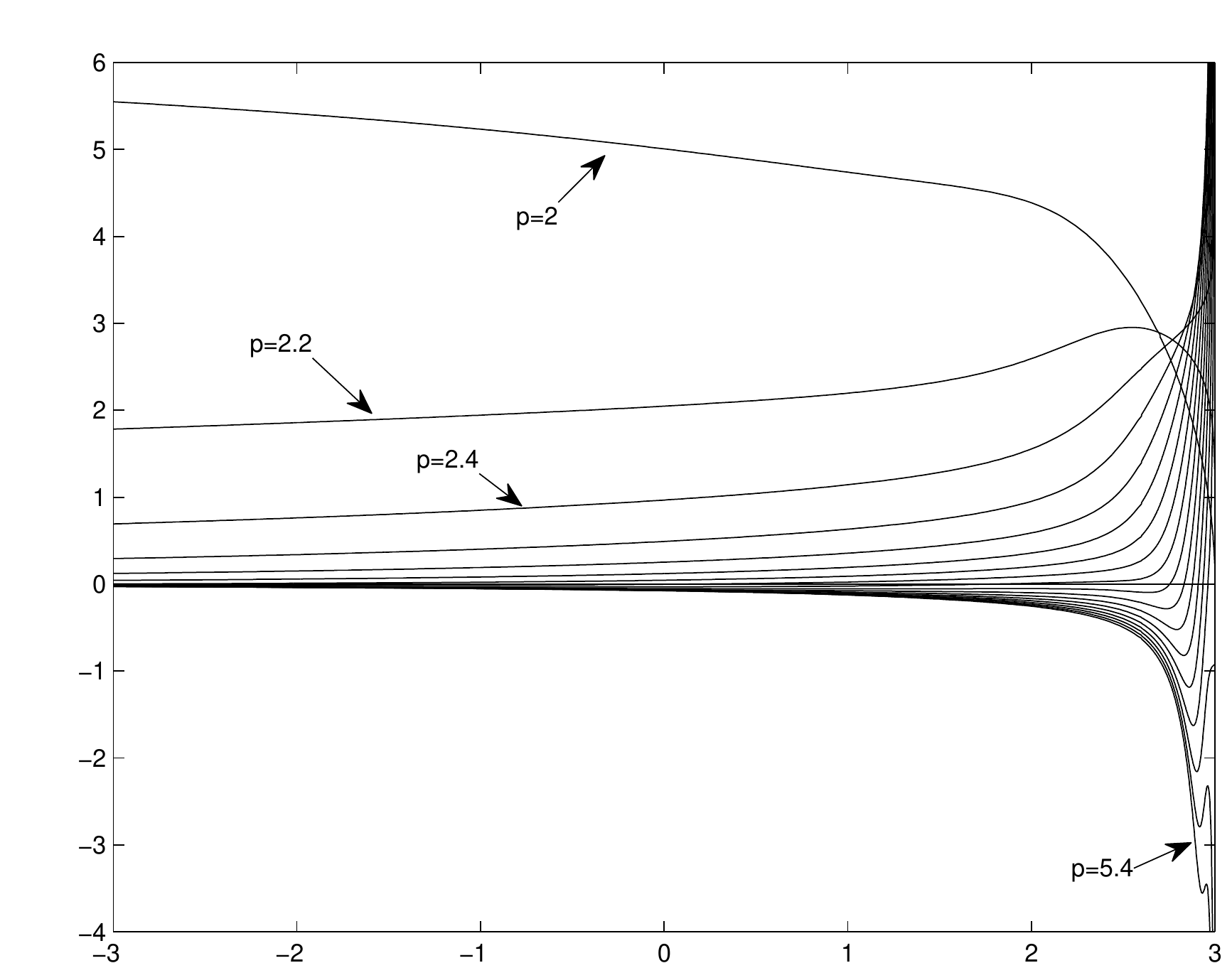}}
   \vskip 5pt
   \scalebox{0.75}{\includegraphics{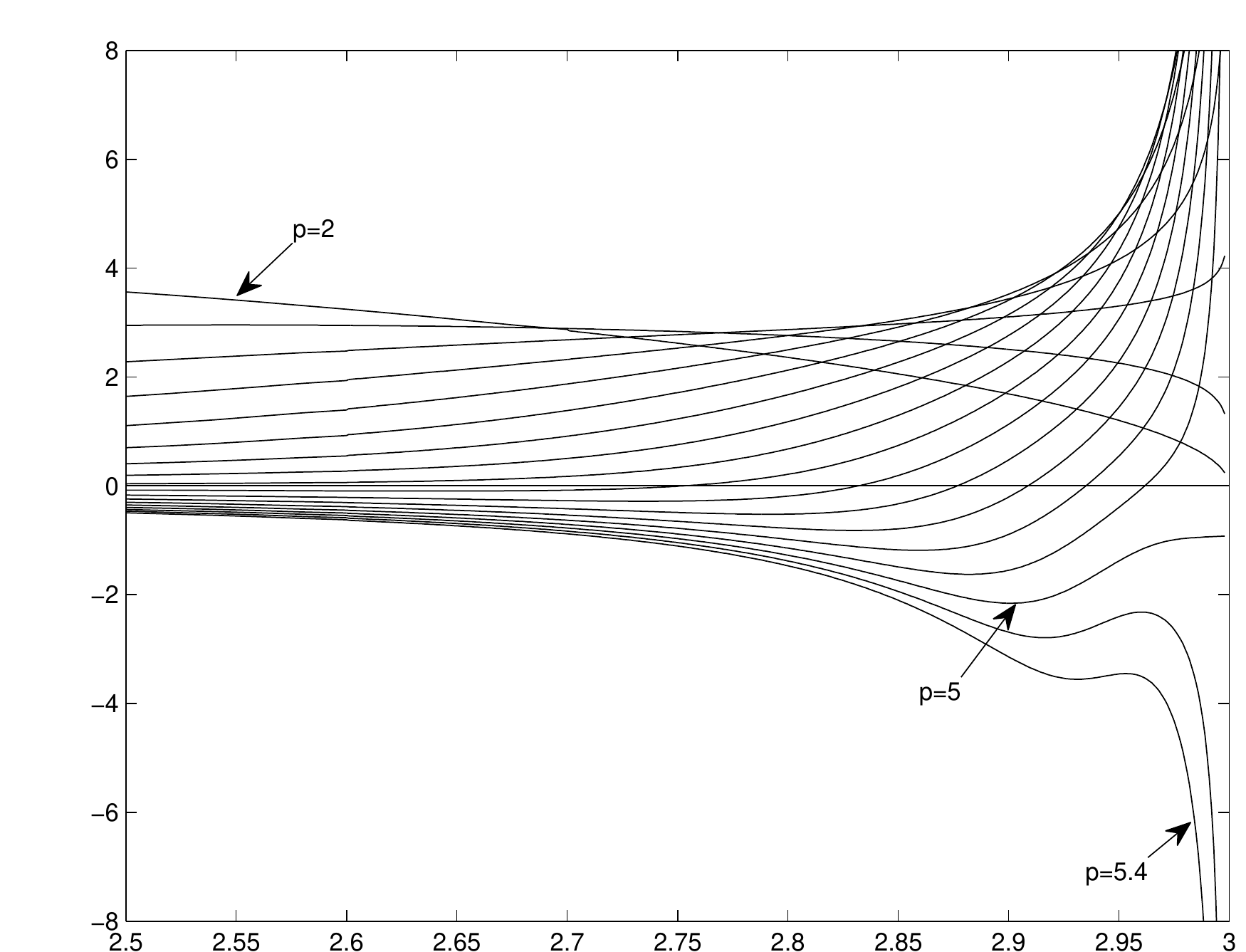}}
  \end{center}
  \caption{Plots of $d_{cc}$ for $f(u)=|u|^{p-1}u$ with with $\beta=2$, $\gamma=1$, $-3\leq c<3$ and $p=2,2.2,2.4,\ldots,5.4$.}\label{F:d1_even}
\end{figure}

\begin{figure}
  \begin{center}
    \scalebox{0.45}{\includegraphics{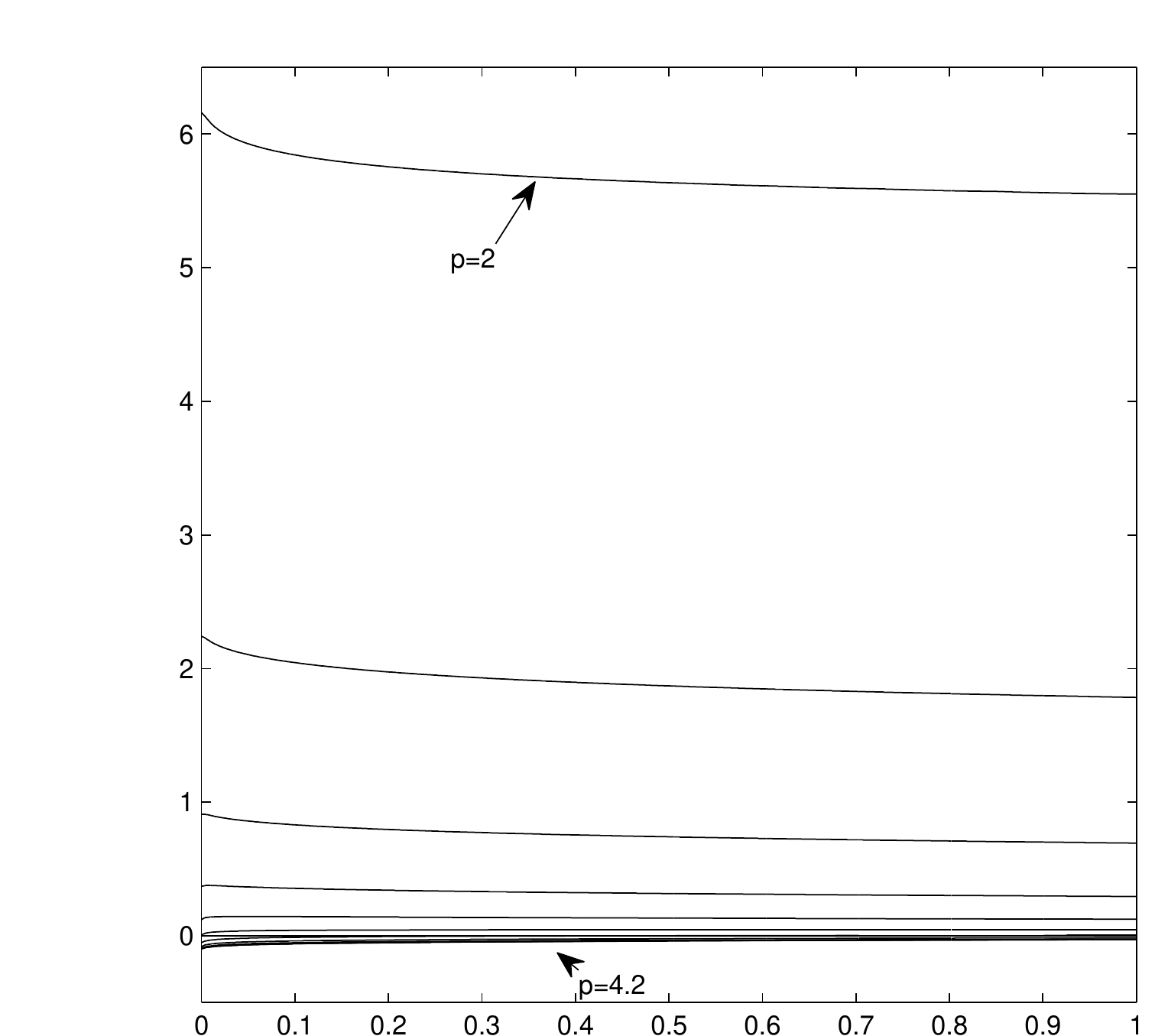}}
    \quad
    \scalebox{0.45}{\includegraphics{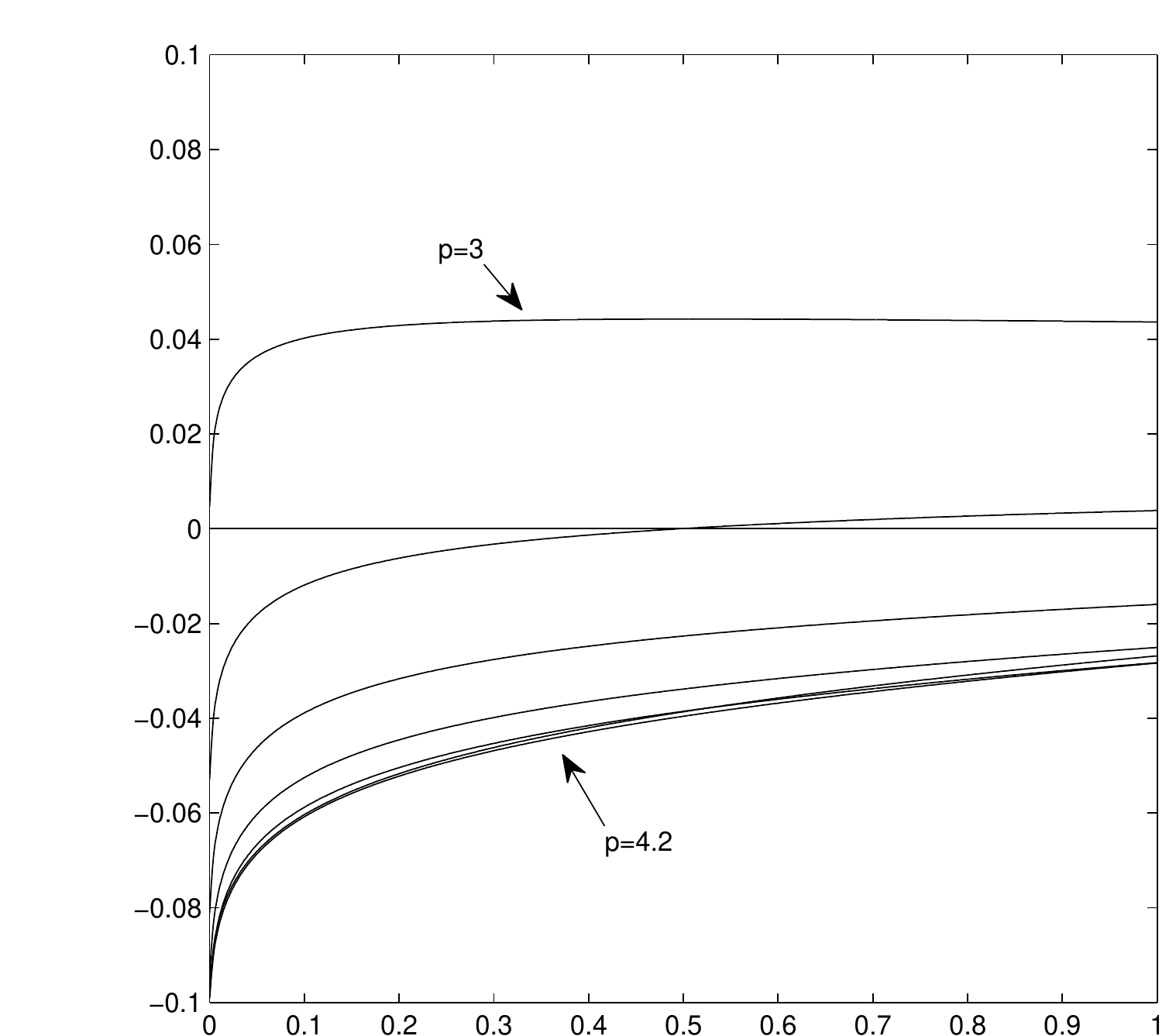}}
  \end{center}
  \caption{Plots of $d_{cc}$ for $f(u)=|u|^{p-1}u$ with $\beta=2$, $c=-3$, $0<\gamma\leq 1$ and $p=2,2.2,2.4,\ldots,4.2$.  The second plot is a blowup of the first, and better illustrates the plots for $3\leq p\leq 4$.}\label{F:d2_even}
\end{figure}

The results for the odd nonlinearity $f(u)=-|u|^{p-1}u$ are shown in Figures \ref{F:d1_even} and \ref{F:d2_even} and summarized in Table \ref{T:d_even}. When $p\leq3$ we have $d_{cc}>0$ for all $c<c_\ast$. On the other hand, when $p\geq 5$ we have $d_{cc}<0$ for all $c< c\ast$. When $3<p<5$ it appears that there exists some speed $c_p$ such that $d_{cc}<0$ for $c<c_p$ and $d_{cc}>0$ for $c_p<c<c_\ast$. Once again, the behavior for small $\gamma>0$ agrees with the results of Theorems \ref{main-instability} and \ref{stability-theo-2}.  The behavior for $c$ near $c_\ast$ is similar to that of the even nonlinearity, only the critical exponent appears to be $p=5$ in this case, in agreement with Theorem \ref{T:d_near_c_star_odd}.

\begin{table}[h!b!p!]
\caption{Sign of $d_{cc}$ for $f(u)=-|u|^{p-1}u$.}
\begin{center}
\begin{tabular}{|l|c|c|}
\hline
  $p$ & Regions where $d_{cc}>0$.\\
  \hline\hline
  2 & $c<c_\ast$ \\
  \hline
  2.2 & $c<c_\ast$ \\
  \hline
  2.4 & $c<c_\ast$\\
  \hline
  2.6 & $c<c_\ast$\\
  \hline
  2.8 & $c<c_\ast$\\
  \hline
  3 & $c<c_\ast$ \\
  \hline
  3.2 & $-1.262c_\ast<c<c_\ast$ \\
  \hline
  3.4 & $0.033c_\ast<c<c_\ast$ \\
  \hline
  3.6 & $0.589c_\ast<c<c_\ast$ \\
  \hline
  3.8 & $0.918c_\ast<c<c_\ast$ \\
  \hline
  4 & $0.944c_\ast<c<c_\ast$ \\
  \hline
  4.2 & $0.959c_\ast<c<c_\ast$ \\
  \hline
  4.4 & $0.970c_\ast<c<c_\ast$ \\
  \hline
  4.6 & $0.978c_\ast<c<c_\ast$ \\
  \hline
  4.8 & $0.987c_\ast<c<c_\ast$ \\
  \hline
  5 & empty \\
  \hline
  \end{tabular}\label{T:d_even}
  \end{center}
\end{table}

\newpage

\newpage

\section*{}
\pdfbookmark[0]{References}{titr-1}

\end{document}